\documentclass[journal]{IEEEtran}
\usepackage{multirow}
\usepackage{xcolor}
\usepackage{ifpdf}
\usepackage{cite}
\usepackage{graphicx}
\usepackage{amsmath}
\usepackage{array}
\usepackage{amssymb}
\usepackage{bmpsize}
\usepackage{algpseudocode,algorithm}
\usepackage{cite}
\usepackage{color}
\usepackage{comment}
\usepackage{mathtools}
\usepackage{booktabs}

\algnewcommand\algorithmicinput{\textbf{Input:}} 
\algnewcommand\Input{\item[\algorithmicinput]}
\algnewcommand\algorithmicoutput{\textbf{Output:}}
\algnewcommand\Output{\item[\algorithmicoutput]}

\newtheorem{remark}{Remark}

\begin{document}
\title{Numerical Spectrum Linking: Identification of Governing PDE via Koopman-Chebyshev Approximation with Resampling}%
\author{
        Phonepaserth~Sisaykeo,~\IEEEmembership{Graduate Student Member,~IEEE,}
        and~Shogo~Muramatsu,~\IEEEmembership{Senior Member,~IEEE}

\thanks{
This work was supported by JSPS KAKENHI Grant Numbers JP22H00512, JP24H00365, JP24K21314 and JP26H02490. 

Phonepaserth Sisaykeo is with the Graduate School of Science and Technology, Niigata University, Niigata 950-2181, Japan.  

Shogo Muramatsu is with the Faculty of Engineering, Niigata University, Niigata 950-2181, Japan (email: shogo@eng.niigata-u.ac.jp). 

Corresponding author: Shogo Muramatsu.

This article has supplementary downloadable material available at [https://codeocean.com/capsule/1315472/tree], provided by the authors. The material includes MATLAB code and experimental configurations used to reproduce the figures and tables presented in this paper. 
}
}
\maketitle

\begin{abstract}
A numerical framework is proposed for identifying governing partial differential equations (PDEs) from observational data by establishing a link between observation-driven and equation-driven Koopman operators in a common Chebyshev spectral domain. In contrast to data-driven approaches such as dynamic mode decomposition (DMD), which approximate Koopman operators without explicitly relating them to differential operators, the proposed framework constructs finite-dimensional Koopman operators using Chebyshev spectral representations, thereby enabling direct comparison between data-derived dynamics and candidate governing PDEs.
A unified observation model together with a least-squares coefficient recovery formulation is introduced to recover Chebyshev spectral coefficients from observations obtained on arbitrary sampling grids. This provides a numerically consistent interface between practical observations and Chebyshev-based Koopman analysis.
Numerical experiments under direct Chebyshev, uniform, and irregular sampling configurations demonstrate that the proposed framework accurately identifies the governing PDE from observations. An observation-density study shows that reliable PDE identification is consistently achieved once sufficient independent observations are available for stable coefficient recovery, providing a practical guideline.
\end{abstract}


\section{Introduction}
\label{sec:intro}

\begin{figure}[tb]
  \centering
  \includegraphics[width=.98\linewidth]{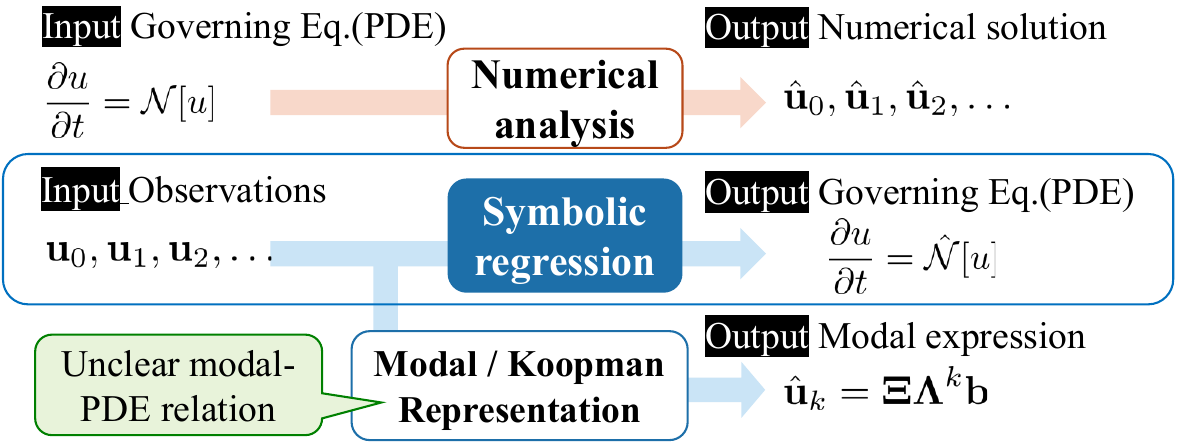}
  \caption{Conceptual motivation for linking different representations of dynamical systems. Numerical analysis, symbolic regression, and modal (Koopman) approaches yield distinct representations, while the relationship between modal representations and governing differential operators remains unclear, motivating the need for a unified framework.}
  \label{fig:background}
\end{figure}

\IEEEPARstart{I}{dentifying} governing PDEs of dynamical systems from observational data is a fundamental problem in signal processing and system identification. In many practical scenarios, system behavior is observed through discrete measurements rather than explicit analytical models, making the recovery of underlying dynamics from sampled data essential for analysis, prediction, and control. While Shannon's sampling theorem provides a theoretical foundation for signal reconstruction from discrete observations~\cite{Shannon1949a,Unser2000a,Eldar2014SamplingSystems}, recent studies have extended the classical sampling framework through Koopman operator theory to characterize broader classes of dynamical signals~\cite{Zeng2024KoopmanSampling}. Extending these sampling concepts to the identification of governing spatiotemporal dynamics remains a significantly more challenging problem.

In many scientific and engineering applications, such as environmental monitoring and fluid dynamics, accurate modeling of system dynamics is essential for prediction and control. For instance, understanding river flow behavior and its interaction with surrounding structures is critical for disaster prevention and environmental management~\cite{Moteki2022CaptureBars,Moteki2023OnSandbars,Ohara2024Physics-informedModel}. Furthermore, new discoveries have overturned established laws~\cite{Karisawa2025Acceleration-inducedSlopes}, emphasizing the importance of data-driven identification of governing dynamics for newly observed phenomena. However, real-world systems are often complex, nonlinear, and partially observed, making it difficult to derive governing PDEs using conventional analytical approaches. This has led to increasing interest in data-driven methodologies for discovering governing equations from measurements.

Data-driven modeling has been advanced by mode-expansion methods such as dynamic mode decomposition (DMD)~\cite{Brunton2022Data-DrivenControl,Bramburger2024Data-DrivenSystems,Strang2019LinearData} and its variants~\cite{Baddoo2023Physics-informedDecomposition,Kobayashi2023Multi-ResolutionModeling}, which approximate the Koopman operator~\cite{2020TheControl,Williams2015ADecomposition,Korda2020Data-drivenOperator}. More recently, learning-based approaches have also been developed to construct Koopman operators directly from observational data using deep neural networks~\cite{Hu2026HybridKoopman}. While these approaches are effective for capturing temporal evolution and constructing surrogate models, the resulting representations are primarily expressed in modal form and lack an explicit connection to differential operators. This limits their interpretability and restricts their applicability for identifying governing PDEs. Therefore, bridging mode expansion with PDE identification, and ultimately symbolic regression~\cite{Wu2025DiscoveringFramework}, is of significant importance, as illustrated in Fig.~\ref{fig:background}.

Recent advances in data-driven differential equation discovery have been comprehensively reviewed by Lou \textit{et al.}~\cite{Lou2026Data-drivenSystems}, highlighting the rapid development of methodologies for identifying governing equations directly from observational data across a wide range of physical systems.
Pioneering work in PDE identification includes PDE-FIND by Rudy \textit{et al.}~\cite{Rudy2017Data-drivenEquations}, which employs sparse regression on candidate libraries of differential operators to successfully identify governing equations of complex systems. However, such approaches rely on predefined libraries, are sensitive to noise, and require reanalysis of the observed data whenever the library is updated. Chen \textit{et al.}'s SGA-PDE~\cite{Chen2022SymbolicSGA-PDE} alleviates library dependence using genetic algorithms, but introduces increased computational cost and stability challenges. More recently, Du \textit{et al.}'s DISCOVER~\cite{Du2024DISCOVER:Learning} utilizes reinforcement learning to identify PDEs with unknown terms, yet issues related to computational efficiency, reward design, and robustness remain. Active-learning-based sparse model discovery has also been investigated to improve data efficiency in the low-data regime~\cite{Larranaga2026HowLow}. On the other hand, physics-informed neural networks (PINNs)~\cite{Raissi2019Physics-informedEquations,Huang2025PartialSurvey} incorporate physical constraints into learning frameworks, but typically assume that governing equations are known a priori and therefore address a different problem setting.
Despite these advances, a fundamental gap remains in establishing a direct and interpretable connection between the Koopman operator and differential operators governing system dynamics. Conventional DMD lacks a clear relationship between the observable basis functions \(\{\varphi_m\}\) and the differential operator \(\mathcal{N}\), limiting its capability for PDE identification from data.

The objective of this study is to establish a unified numerical framework that bridges this gap. Motivated by the use of Chebyshev polynomial bases in Koopman-based representations~\cite{Shi2024KoopmanStudyb}, a framework termed \emph{numerical spectrum linking} is developed. Chebyshev polynomials approximate observed functions through orthogonal basis expansions, while discrete cosine transform (DCT) enables efficient finite-dimensional representations~\cite{Ahmed1974DiscreteTransform,Ochoa-Dominguez2019DiscreteTransform}. This allows differential operators to be represented as matrix operations on coefficient vectors~\cite{Barrio2004AlgorithmsSeries}. Meanwhile, Koopman mode expansion provides a linear representation of nonlinear dynamical systems. By connecting these two representations, the proposed framework enables identification of governing PDEs through spectral analysis.

This paper is an extension of a preliminary conference version presented in~\cite{Sisaykeo2026NumSpecLink}, where the numerical spectrum linking framework was formulated under the ideal assumption that observations are directly available at Chebyshev sampling nodes. 
However, this assumption constitutes a restrictive condition in practical measurement settings. In practical applications, observations are typically acquired on uniform grids or irregular spatial locations, leading to a mismatch between the sampling domain and the spectral representation required for Chebyshev-based analysis. This sampling mismatch degrades coefficient estimation accuracy and consequently affects Koopman operator estimation and PDE identification performance.

To address this issue, the present study introduces a unified observation model and coefficient recovery framework that enables observations acquired on uniform or irregular sampling grids to be consistently represented in the Chebyshev spectral domain. This generalized formulation enables consistent Chebyshev spectral analysis and Koopman operator estimation under non-ideal observation conditions. In addition, expanded analyses are provided to investigate the influence of sampling mismatch and coefficient recovery errors on PDE identification performance.

The preliminary conference version relied on eigenvalue-based and eigenvector-based comparison metrics to assess agreement between the observation-driven and equation-driven Koopman operators. In this extended study, we identify a structural property of the equation-driven operator: its generator matrix is nilpotent for every candidate PDE considered, forcing a degenerate unit eigenvalue independent of the underlying dynamics, which renders such eigenvalue-based comparisons uninformative for PDE identification. This motivates the data-projected discrepancy criterion adopted throughout this paper, which instead compares the action of the two operators on observed data.

The main aspects of the proposed study are summarized as follows.
\begin{itemize}
    \item First, a unified numerical framework is established for linking Chebyshev spectral representations with Koopman operator theory for PDE identification through a common spectral-domain formulation.
    \item Second, a unified observation model and coefficient recovery formulation is incorporated to enable consistent Chebyshev spectral analysis from observations acquired on uniform and irregular sampling grids.
    \item Third, the influence of sampling mismatch and resampling on Koopman operator estimation and PDE identification accuracy is systematically evaluated under multiple observation settings.
    \item Fourth, a structural nilpotency property of the equation-driven generator matrix is identified, which forces a degenerate eigenvalue independent of the candidate PDE and renders purely eigenvalue-based comparison uninformative, motivating the data-projected discrepancy criterion adopted for PDE identification throughout this study.
\end{itemize}

The remainder of this paper is organized as follows. Section~\ref{sec:theories} reviews the foundational theories of Koopman operators and Chebyshev approximation. Section~\ref{sec:proposed} presents the proposed numerical spectrum linking framework with resampling. Section~\ref{sec:performance} provides numerical experiments under various sampling conditions. Finally, Section~\ref{sec:conclusion} concludes the paper.

\section{Overview of Foundational Theories}
\label{sec:theories}

\begin{figure*}[tb]
  \centering
  \includegraphics[width=.75\linewidth]{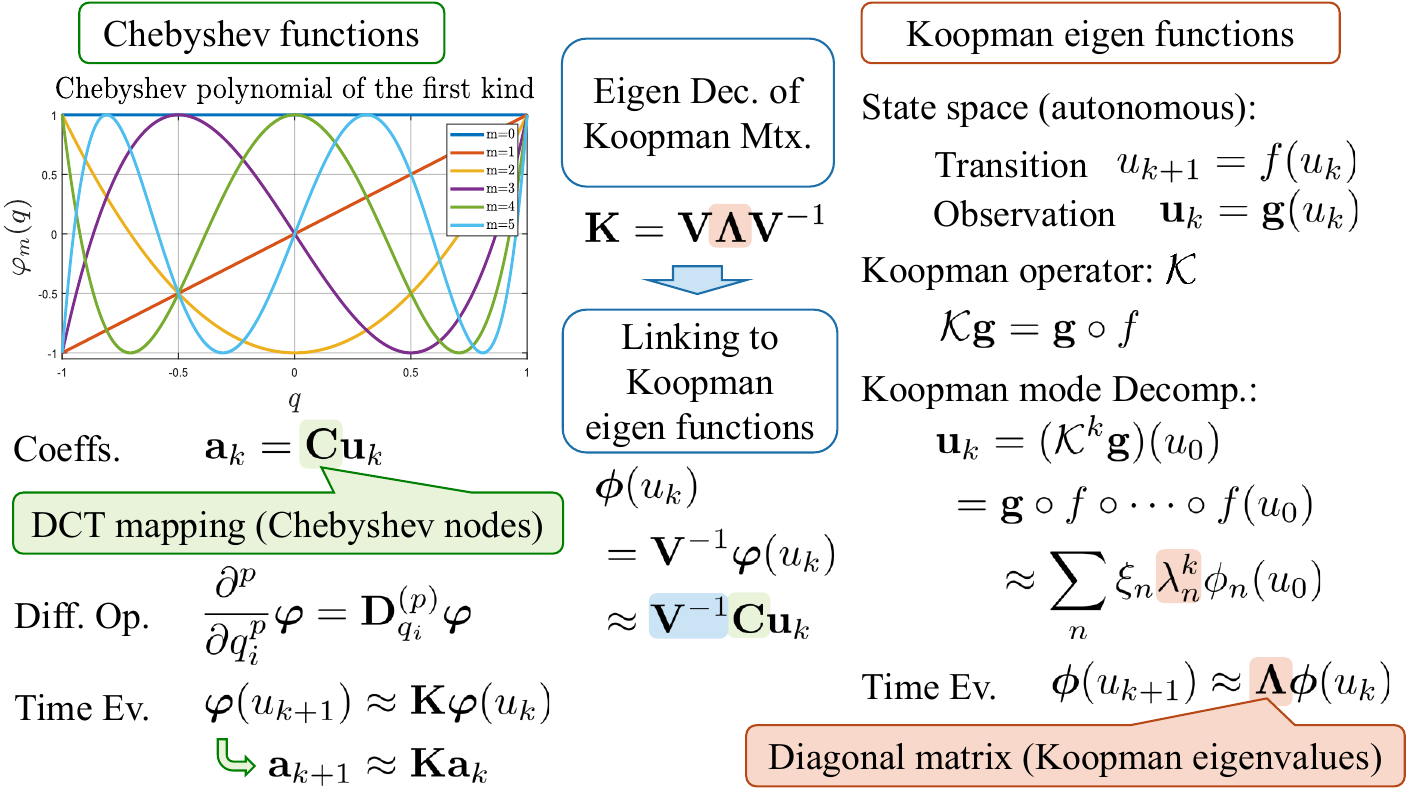}
  \caption{Relationship between Chebyshev spectral representation and Koopman operator theory. The system state is mapped to Chebyshev coefficients via DCT, while Koopman operator theory provides a modal decomposition of the dynamics. Their connection is established through spectral decomposition of the Koopman operator, forming the theoretical basis of the proposed framework.}
  \label{fig:outline}
\end{figure*}

In this section, we review the foundational theories underlying the proposed framework. 
We first introduce the Koopman operator in discrete time and its connection to data-driven modal analysis, which provides a linear representation of nonlinear dynamics. 
We then present Chebyshev polynomials and their relation to  DCT, which form the basis of efficient spectral approximation. 
Finally, we summarize differentiation in the Chebyshev basis, through which differential operators can be represented as matrix operations on coefficient vectors. 
Together, these concepts establish the mathematical foundation for constructing and comparing observation-driven and equation-driven operators, as illustrated in Fig.~\ref{fig:outline}.

\subsection{Koopman Operator}

The description of dynamical systems often benefits from a state-space representation, which is compatible with modern control theory. 
In continuous time, such systems are described by ordinary differential equations (ODEs) or PDEs, while in discrete time they are described by difference equations or maps. 
These representations are directly connected to optimal control and state estimation.

A theoretical framework for analyzing such systems is provided by the \emph{Koopman operator}~\cite{2020TheControl,Strasser2025AnGuarantees}. 
The Koopman operator enables one to treat a nonlinear dynamical system as a linear operator acting on observables. 
Consider the nonlinear discrete-time autonomous system
\begin{subequations}
\begin{align}
u_{k+1} &= f(u_k), \quad u_k\in\mathcal{M},\\
{y}_{k} &= {g}(u_{k}), \quad {y}_{k}\in\mathbb{R},
\end{align}
\end{subequations}
where $k\in\mathbb{N}_0$ denotes discrete time, $\mathcal{M}$ is the state space, 
$f:\mathcal{M}\to\mathcal{M}$ is a nonlinear (regular) state-transition map, 
and ${g}:\mathcal{M}\to\mathbb{R}$ is an observable.

Let $\mathcal{M}^\ast$ denote the space of admissible observables ${g}$. 
The Koopman operator $\mathcal{K}:\mathcal{M}^\ast\to\mathcal{M}^\ast$ is defined by composition as
\begin{equation}
  \mathcal{K}{g} \coloneqq {g}\circ f, \quad \forall g\in\mathcal{M}^\ast.    
\end{equation}
Thus, instead of describing the evolution of the state directly, the Koopman operator describes the evolution of observables. 
Note that $\mathcal{K}$ is linear, i.e.,
\(
\mathcal{K}(a {g}_1 + b {g}_2) = a \mathcal{K}{g}_1 + b \mathcal{K}{g}_2,
\)
even though the underlying dynamics $f$ is nonlinear. 
Unless $\mathcal{M}$ is finite, however, the operator $\mathcal{K}$ is infinite-dimensional. 
The operator can also be defined for vector-valued observables $\mathbf{g}:\mathcal{M}\to\mathbb{R}^M$.

\subsubsection{Koopman Mode Expansion}

If the system is integrable, the state space is compact, and scalar observables are considered, 
the spectral properties of $\mathcal{K}$ allow for an eigenfunction expansion. 
A function $\phi_\lambda\in\mathcal{M}^\ast\setminus\{0\}$ satisfying
\begin{equation}
  \mathcal{K}\phi_\lambda = \phi_\lambda\circ f = \lambda \phi_\lambda
\end{equation}
is called a \emph{Koopman eigenfunction} with corresponding \emph{Koopman eigenvalue} $\lambda\in\mathbb{C}$. 
Using a family of eigenfunctions $\{\phi_{\lambda_m}\}_{m=0}^\infty$, an observable $g$ can be expanded as
\(
g = \sum_{m=0}^\infty \xi_m \phi_{\lambda_m},
\)
where $\xi_m\in\mathbb{C}$ are the \emph{Koopman modes}. 
The time evolution of $g$ along a trajectory is then
\(
y_{k} = (\mathcal{K}^{k}g)(u_0)
= \sum_{m=0}^\infty \xi_m \lambda_m^{k} \phi_{\lambda_m}(u_0).
\)

For vector-valued observables $\mathbf{g}$, both the measurements 
$\mathbf{y}_{k}\in\mathbb{R}^M$ 
and the Koopman modes $\boldsymbol{\xi}_m\in\mathbb{C}^M$ become vector-valued, yielding
\(
\mathbf{y}_{k} = (\mathcal{K}^{k}\mathbf{g})(u_0) 
= \sum_{m=0}^\infty \boldsymbol{\xi}_m \lambda_m^{k} \phi_{\lambda_m}(u_0).
\)

This expansion provides a linear spectral description of nonlinear dynamics and motivates the finite-dimensional operator representations used later.

\subsubsection{Koopman Matrix Approximation}

For numerical computation, the infinite-dimensional operator $\mathcal{K}$ is approximated 
in a finite-dimensional subspace spanned by a set of basis functions 
$\{\varphi_m\}_{m=0}^{M-1}\subseteq\mathcal{M}^\ast$, which may not be a set of the Koopman eigenfunctions. 
Approximating $g(u_k)$ by a linear combination of these basis functions leads to the relation
\begin{equation}
(\mathcal{K}g)(u_k) \approx \sum_{m=0}^{M-1} a_m (\mathcal{K}\varphi_m)(u_k)
= \sum_{m=0}^{M-1} a_m \varphi_m(u_{k+1}). \label{eq:koopmanmatrix0}
\end{equation}
If each $\varphi_m(u_{k+1})$ can be expressed as a linear combination of $\{\varphi_n(u_k)\}$, 
the system admits a matrix representation
\begin{equation}
\boldsymbol{\varphi}(u_{k+1}) \approx \mathbf{K}\,\boldsymbol{\varphi}(u_{k}),
\end{equation}
where $\boldsymbol{\varphi}(u) = (\varphi_0(u),\ldots,\varphi_{M-1}(u))^\intercal$ and
$\mathbf{K}\in\mathbb{C}^{M\times M}$ is called the \emph{Koopman matrix}. 
Given data $\{u_{k}\}_{k=0}^{N-1}$,
$\mathbf{K}$ can be estimated by solving the least-squares problem
\begin{equation}
\hat{\mathbf{K}} \in \arg\min_{\mathbf{K}} \sum_{k=0}^{N-2} 
\left\| \boldsymbol{\varphi}(u_{k+1}) - \mathbf{K}\boldsymbol{\varphi}(u_{k})\right\|_2^2.
\end{equation}
If the basis functions are Koopman eigenfunctions, the matrix $\mathbf{K}$ becomes diagonal. 
Eigenvalue analysis of the estimated $\mathbf{K}$ provides access to Koopman modes, eigenvalues, and eigenfunctions. 
This finite-dimensional approximation procedure underlies DMD and its extensions~\cite{Brunton2022Data-DrivenControl,Bramburger2024Data-DrivenSystems,Strang2019LinearData,Bistrian2025ReducedLearning}. 
In the proposed framework, such finite-dimensional Koopman matrices will be constructed both from observed data and from equation-driven spectral representations. Chebyshev polynomials provide the link between the two.

\subsection{Chebyshev Polynomials}

In this subsection, we summarize the basic properties of Chebyshev polynomials and their role in spectral approximation. 
Chebyshev polynomials form an orthogonal basis on $[-1,1]$ with respect to a specific weight function, and their tensor-product extension enables efficient approximation of multivariate functions~\cite{Mason2002ChebyshevPolynomials}. 
By sampling at Chebyshev nodes and exploiting the correspondence with Type-II DCT (DCT-II), the expansion coefficients can be computed efficiently~\cite{Ochoa-Dominguez2019DiscreteTransform}. 
This connection provides a practical tool for constructing finite-dimensional operator representations and serves as a foundation for the Koopman operator analysis developed later.

\subsubsection{Multidimensional Chebyshev basis and expansion}

We consider a real-valued function on a $D$-dimensional hypercube,
\begin{equation}
u \colon [-1,1]^D \to \mathbb{R} \colon \mathbf{q} \mapsto u(\mathbf{q}).
\end{equation}
For simplicity, the time index is omitted.

Let $T_m(q)\coloneqq \cos(m\arccos q)$ be the Chebyshev polynomials of the first kind, $m\in\mathbb{N}_0$, $q\in[-1,1]$. 
We define the $D$-dimensional tensor-product basis
\(
\tau_{\mathbf{m}}(\mathbf{q}) = \prod_{d=1}^D T_{m_d}(q_d),
\)
where
\(
\mathbf{m}=(m_1,\dots,m_D)^\intercal\in\mathbb{N}_0^D.
\)
Then, $u$ admits the formal Chebyshev expansion
\(
u(\mathbf{q}) = \sum_{\mathbf{m}\in\mathbb{N}_0^D} a_\mathbf{m}\,\tau_{\mathbf{m}}(\mathbf{q}).
\)

This tensor-product construction provides a natural spectral basis for multidimensional state representations and is particularly suitable for PDE-related dynamics on bounded domains~\cite{Karunakar2019ShiftedEquations}.

\subsubsection{Sampling at Chebyshev nodes and DCT-II computation}

For a finite-dimensional approximation, we define 
\(
\mathcal{D}(\mathbf{M}) \coloneq \{0,1,\dots,M_1-1\}
\times\cdots\times
\{0,1,\dots,M_D-1\}
\)
with
\(
\mathbf{M}=\mathrm{diag}(M_1,M_2,\dots,M_D),
\)
and write
\begin{equation}
u(\mathbf{q}) \approx \sum_{\mathbf{m}\in\mathcal{D}(\mathbf{M})} a_\mathbf{m}\,\tau_{\mathbf{m}}(\mathbf{q})
= \sum_{\mathbf{m}\in\mathcal{D}(\mathbf{M})} \check{a}_\mathbf{m}\,\check{\tau}_{\mathbf{m}}(\mathbf{q}).
\label{eq:cheb_finite_expansion}
\end{equation}
Here, we rescale with
\(
\check{a}_\mathbf{m} = \gamma_{\mathbf{m}}^{-1} a_\mathbf{m}, \quad
\check{\tau}_{\mathbf{m}}(\cdot) = \gamma_{\mathbf{m}} \tau_{\mathbf{m}}(\cdot), \quad
\gamma_{\mathbf{m}}=\prod_{d=1}^D \gamma^{(d)}_{m_d},
\)
and choose per-dimension factors 
$\gamma^{(d)}_0 = M_d^{-1/2}$ and $\gamma^{(d)}_{m\ge 1} = (2/M_d)^{1/2}$,
so that the discrete transform below becomes orthonormal. 
With this normalization, the coefficients are expressed as
\(
\check{a}_\mathbf{m}
=
{\langle \check{\tau}_{\mathbf{m}}, u \rangle_w}/{\langle \check{\tau}_{\mathbf{m}}, \check{\tau}_{\mathbf{m}} \rangle_w}
=
(\prod_{d=1}^D \frac{M_d}{\pi})\,\langle \check{\tau}_{\mathbf{m}}, u \rangle_w
\),
where $w(\mathbf{q})=\prod_{d=1}^{D}(1-q_d^2)^{-1/2}$ denotes the Chebyshev weight~\cite{Mason2002ChebyshevPolynomials,Ochoa-Dominguez2019DiscreteTransform}.

Assume each $M_d$ is even and use interior Chebyshev nodes (without endpoints)
\begin{equation}
p_{n_d} \coloneqq \cos\!\left(\frac{(2n_d+1)\pi}{2M_d}\right), 
\quad n_d=0,\dots,M_d-1.
\label{eq:cheb_interior_nodes}
\end{equation}
Let
\begin{equation}
\mathbf{p}_{\mathbf{n}} \coloneqq (p_{n_1},\dots,p_{n_D})^\intercal,
\label{eq:cheb_node_vector}
\end{equation}
and sample vector
\(
\mathbf{u} \coloneqq \bigl(u(\mathbf{p}_{\mathbf{n}})\bigr)_{\mathbf{n}\in\mathcal{D}(\mathbf{M})}.
\)
Then, the coefficient vector 
\(
\mathbf{a}=\bigl(\check{a}_\mathbf{m}\bigr)_{\mathbf{m}\in\mathcal{D}(\mathbf{M})}
\)
is well approximated by a $D$-dimensional DCT-II:
\begin{equation}
\mathbf{a} \approx \mathbf{C}\mathbf{u},\label{eq:a_approx_Cu}
\end{equation}
where $\mathbf{C}$ denotes the Kronecker-structured matrix of the $D$-dimensional DCT-II. 
We can define the analysis map as
\(
\boldsymbol{\varphi}(u) \approx \mathbf{C}\mathbf{u},
\)
which sends samples on Chebyshev nodes to scaled Chebyshev coefficients.
Let
\(
u_{k}(\mathbf{q})
\coloneqq
u(\mathbf{q},k\Delta t)
\)
denote the system state at the discrete time
\(t=k\Delta t\).
The corresponding samples on Chebyshev nodes are collected in
\(
\mathbf{u}_{k}
=
\left(
u_{k}(\mathbf{p}_{\mathbf n})
\right)_{\mathbf n\in\mathcal D(\mathbf M)}.
\)
Accordingly, the observation on the Chebyshev nodes can be expressed as
\begin{equation}
\mathbf{u}_{k}
=
\mathbf{g}(u_{k})
\approx
\mathbf{C}^{\intercal}\boldsymbol{\varphi}(u_{k}).
\end{equation}

This DCT-based coefficient computation is essential in the proposed framework because it provides the spectral coordinates in which both observation-driven and equation-driven operators are represented and compared.

\subsection{Differentiation of Chebyshev Polynomials}

We summarize differentiation identities and their implementation in both coefficient space and node space, consistent with Chebyshev polynomial differentiation theory~\cite{Prodinger2017RepresentingQuestions}.

For the $D$-dimensional tensor-product basis $\tau_{\mathbf{m}}(\mathbf{q})$ with 
$\mathbf{q}=(q_d)_{d=1}^{D}$, the partial derivative with respect to $q_d$ acts on the coefficient vector as a Kronecker product operator
\[
\mathbf{D}^{(d)} 
= \mathbf{I}_{M_D}\otimes\cdots\otimes \mathbf{I}_{M_{d+1}}
\otimes \mathbf{D}_{M_d}
\otimes \mathbf{I}_{M_{d-1}}\otimes\cdots\otimes \mathbf{I}_{M_1},
\]
so that for the stacked coefficients ${\mathbf{a}}\in\mathbb{R}^{M_1\cdots M_D}$,
\begin{equation}
\frac{\partial }{\partial q_d}u(\mathbf{q})
\approx
\sum_{\mathbf{m}\in\mathcal{D}(\mathbf{M})} \check{a}_\mathbf{m}'\check{\tau}_\mathbf{m}(\mathbf{q})
\quad\Longleftrightarrow\quad 
\mathbf{a}' \approx \mathbf{D}^{(d)}\,\mathbf{a},
\label{eq:mtxrep_diff}
\end{equation}
where $\mathbf{a}\approx\mathbf{C}\mathbf{u}$ and $\mathbf{a}'\approx\mathbf{C} (\partial\mathbf{u}/\partial q_d)$
with $\partial\mathbf{u}/\partial q_d$ denoting the derivative sampled at the same nodes as
$\mathbf{u}$~\cite{Bedratyuk2022DerivationsPolynomials}.
Then, the partial derivative evaluated at nodes is computed by
\begin{equation}
\frac{\partial \mathbf{u}}{\partial q_d}
\coloneqq
\left(\frac{\partial u(\mathbf{p}_{\mathbf{n}})}{\partial q_d}\right)_{\mathbf{n}\in\mathcal{D}(\mathbf{M})}
\approx
\mathbf{C}^\intercal\,\mathbf{D}^{(d)}\,\mathbf{C}\,\mathbf{u}.
\end{equation}

Second and higher derivatives are obtained by repeated application in coefficient space. For example,
\(
\frac{\partial^2 }{\partial q_d \partial q_r}u(\mathbf{q})
\approx
\sum_{\mathbf{m}\in\mathcal{D}(\mathbf{M})} \check{a}_\mathbf{m}''\check{\tau}_\mathbf{m}(\mathbf{q})
\)
is equivalent to
\(
\mathbf{a}''\approx\mathbf{D}^{(d)}\mathbf{D}^{(r)}\mathbf{a},
\)
which in node space corresponds to 
\(
\mathbf{C}^\intercal \mathbf{D}^{(d)}\mathbf{D}^{(r)} \mathbf{C}.
\)

This formulation enables differential operators to be represented in matrix form on the same spectral space used for Koopman approximation. 
Consequently, it provides the bridge needed to derive equation-driven operators from governing PDEs and to compare them with data-driven operators in the proposed numerical spectrum linking framework.

\section{Proposed Numerical Spectrum Linking}
\label{sec:proposed}

\begin{figure}[tb]
  \centering
  \includegraphics[width=.98\linewidth]{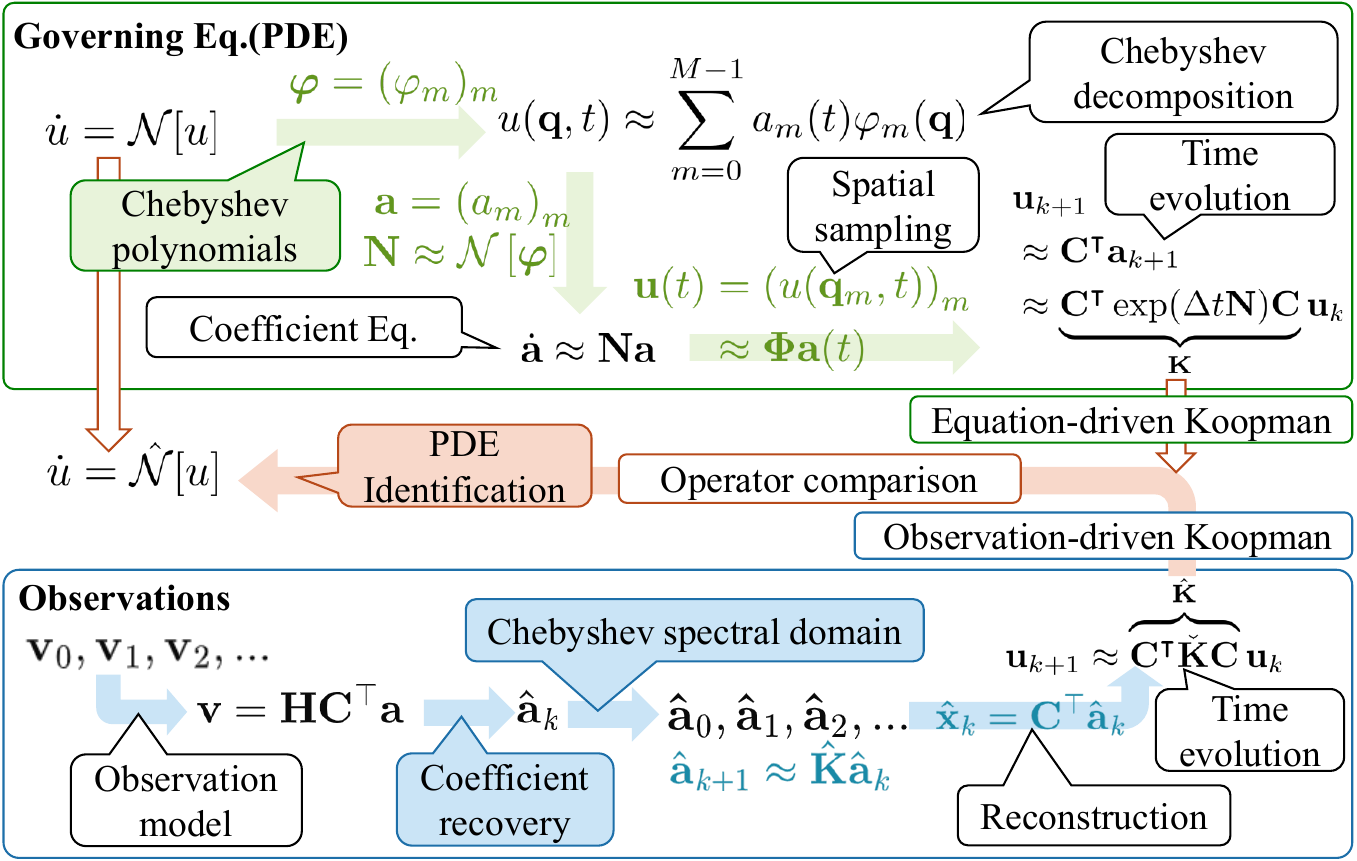}
  \caption{Overview of the proposed numerical spectrum linking framework. Observed data are related to the Chebyshev spectral representation through an observation model, and Chebyshev coefficients are recovered from the observations. From the recovered coefficient sequence, an observation-driven Koopman matrix $\hat{\mathbf{K}}$ is constructed. In parallel, an equation-driven Koopman matrix $\mathbf{K}^{\star}$ is derived from the governing PDE using a Chebyshev-based representation. The governing equation is identified by comparing the two operators.}
  \label{fig:proposal}
\end{figure}

We present the proposed numerical spectrum linking framework for identifying governing dynamics from observational data. 
The key idea is to construct Koopman matrices from both data and governing equations in a common spectral space, and to compare their spectral properties for consistency, as illustrated in Fig.~\ref{fig:proposal}.

Consider a governing equation of the form
\begin{equation}
    \dot{u} = \mathcal{N}[u],\label{eq:governingeq}
\end{equation}
where $\dot{u}\coloneqq \partial u/\partial t$, and $\mathcal{N}[\cdot]$ is a spatial differential operator. 
Using the Chebyshev-based representation introduced in Section~\ref{sec:theories}, both observation data and differential operators can be expressed in a unified spectral domain.

\subsection{Observation-Driven Identification of Koopman Matrix}
\label{subsec:observation-driven}
Consider a sequence of observations
\(
\{\mathbf{v}_k\}_{k=0}^{N-1}
\)
obtained from a dynamical system governed by \eqref{eq:governingeq}, where
\(
\mathbf{v}_k=\mathbf{g}(u_k)
\)
denotes the observation vector associated with the system state $u_k$ at time step $k$.

Since the observations are not necessarily acquired on Chebyshev nodes, the corresponding Chebyshev coefficient vector is first estimated through the observation model and coefficient recovery procedure described in Sections~\ref{subsec:observation-model-under-arbitrary-sampling} to~\ref{subsec:coefficient-recovery}.

Let
\(
\{\hat{\mathbf{a}}_k\}_{k=0}^{N-1}
\)
denote the resulting sequence of recovered Chebyshev coefficient vectors. The observation-driven Koopman matrix is estimated by solving the least-squares problem
\begin{equation}
    \hat{\mathbf{K}}
    \in
    \arg\min_{\mathbf{K}}
    \sum_{k=0}^{N-2}
    \left\|
    \hat{\mathbf{a}}_{k+1}
    -
    \mathbf{K}\hat{\mathbf{a}}_k
    \right\|_2^2.
\end{equation}
Equivalently, the solution is given by
\begin{equation}
\hat{\mathbf{K}}
=
\hat{\mathbf{A}}_1
\hat{\mathbf{A}}_0^{\dagger},
\label{eq:koopman_estimation}
\end{equation}
where
\(
\hat{\mathbf{A}}_0
\coloneqq
\left(
\hat{\mathbf{a}}_0\;
\hat{\mathbf{a}}_1\;
\cdots\;
\hat{\mathbf{a}}_{N-2}
\right)
\)
and
\(
\hat{\mathbf{A}}_1
\coloneqq
\left(
\hat{\mathbf{a}}_1\;
\hat{\mathbf{a}}_2\;
\cdots\;
\hat{\mathbf{a}}_{N-1}
\right)
\)
are snapshot matrices constructed from the recovered coefficient sequence.

This procedure yields a finite-dimensional approximation of the Koopman operator constructed directly from observational data, similar to recent data-driven Koopman-based system identification approaches~\cite{Ketthong2024Data-DrivenDisturbance}. By expressing the dynamics in the Chebyshev spectral space, the estimated operator $\hat{\mathbf{K}}$ can be directly compared with equation-driven operators derived from the governing PDE.

\subsection{Equation-Driven Derivation of Koopman Matrix}
\label{subsec:equation-driven}

In addition to the observation-driven formulation, a Koopman matrix can be derived directly from the governing equation. 
Consider the continuous-time dynamical system described by \eqref{eq:governingeq}. 
Using the Chebyshev spectral representation introduced in Section~\ref{sec:theories}, the function $u(\mathbf{q},t)$ can be approximated by a finite-dimensional expansion with coefficient vector $\mathbf{a}(t)$.

Applying the differential operator $\mathcal{N}[\cdot]$ to this expansion and using the matrix representation of differentiation in \eqref{eq:mtxrep_diff}, the governing equation can be expressed in coefficient space as
\begin{equation}
\frac{\partial u}{\partial t}\approx
\mathcal{N}\left[
\sum_{\mathbf{m}\in\mathcal{D}(\mathbf{M})} \check{a}_\mathbf{m}\check{\tau}_{\mathbf{m}}(\mathbf{q})\right]
\quad\Longleftrightarrow\quad 
\dot{\mathbf{a}}(t)\approx \mathbf{N}\mathbf{a}(t),
\label{eq:coefficient_ode}
\end{equation}
where $\mathbf{N}$ denotes the matrix representation of the spatial differential operator $\mathcal{N}$ with respect to the Chebyshev basis.

This formulation provides a linear ODE in the spectral domain. 
By discretizing in time with step size $\Delta t$, the evolution of the coefficient vector is approximated as
\begin{equation}
  \mathbf{a}_{k+1} \approx \mathbf{K}^{\star}\mathbf{a}_k,
\end{equation}
where the equation-driven Koopman matrix is given by
\begin{equation}
\mathbf{K}^{\star} \coloneqq \exp(\Delta t\, \mathbf{N}).
\end{equation}

Thus, the Koopman matrix can be constructed directly from the governing differential equation through its spectral representation. 
This equation-driven matrix $\mathbf{K}^{\star}$ serves as a reference for the underlying dynamics and will be compared with the observation-driven matrix $\hat{\mathbf{K}}$ through spectral analysis in the subsequent subsection.

\subsection{Numerical Spectrum Linking}
\label{subsec:numerical-spectrum-linking}

From Sections~\ref{subsec:observation-driven} and~\ref{subsec:equation-driven}, we obtain two approximations of the Koopman operator: 
the observation-driven operator, denoted by 
\(
\hat{\mathbf{K}},
\)
and the equation-driven matrix derived from the governing PDE, denoted by 
\(
\mathbf{K}^{\star}.
\)
Although these matrices are constructed from different principles, they are both represented in the same Chebyshev spectral space, enabling a consistent comparison.

A natural starting point is to compare the two operators through their spectral properties. Let the eigenvalue decompositions of the two matrices be given by
\begin{equation}
  \hat{\mathbf{K}}
  \hat{\mathbf{r}}_j 
  = \hat{\lambda}_j
  \hat{\mathbf{r}}_j, 
  \label{eq:koopman_eig_hat}
\end{equation}
\begin{equation}
  \mathbf{K}^{\star} \mathbf{r}_j^{\star} 
  = \lambda_j^{\star} \mathbf{r}_j^{\star},
  \label{eq:koopman_eig_star}
\end{equation}
where $\hat{\lambda}_j$ and $\lambda_j^{\star}$ denote the Koopman eigenvalues, 
and $\hat{\mathbf{r}}_j$ and $\mathbf{r}_j^{\star}$ denote the corresponding eigenvectors. If the two spectra and their associated Koopman modes agreed, this would suggest that $\hat{\mathbf{K}}$ reproduces both the temporal evolution and the spatial organization encoded in $\mathbf{K}^{\star}$. The following structural property shows, however, that such a comparison is uninformative for the candidate PDEs considered in this study.

\begin{figure}[tb]
  \centering
  \includegraphics[trim=15pt 5pt 20pt 5pt, clip, width=0.98\linewidth]{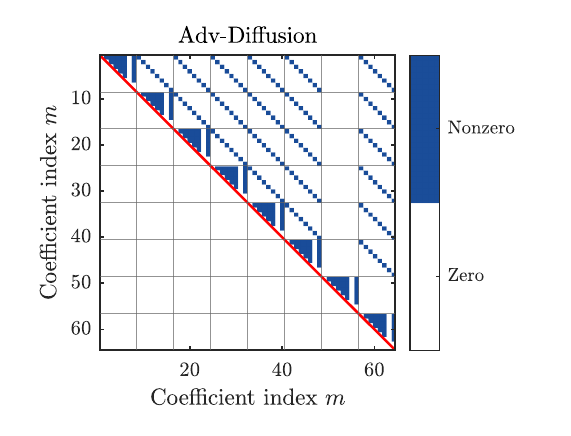}
  \caption{Nonzero pattern of the generator matrix $\mathbf{N}$ for the Advection-Diffusion candidate. All nonzero entries (blue) lie above the main diagonal (red line), confirming the upper triangular, zero-diagonal structure underlying Remark~\ref{rem:nilpotent}.}
  \label{fig:generator_structure}
\end{figure}

\begin{remark}
\label{rem:nilpotent}
For every candidate governing operator $\mathcal{N}[\cdot]$ considered in this study (advection, diffusion, and their combination), the associated coefficient-space matrix $\mathbf{N}$ in \eqref{eq:coefficient_ode} is strictly upper triangular, regardless of the specific coefficient values (e.g., advection speed or diffusivity) entering $\mathcal{N}$. This follows directly from the structure of Chebyshev differentiation: differentiating the $m$-th basis polynomial strictly lowers its degree, so the per-dimension differentiation matrix $\mathbf{D}^{(d)}$ has nonzero entries only strictly above the diagonal. Since $\mathbf{N}$ is built from Kronecker products and linear combinations of such matrices with identity matrices, it inherits the same strictly upper triangular, zero-diagonal structure with respect to the coefficient ordering (Fig.~\ref{fig:generator_structure}). This argument extends to derivatives of any order: since a product of strictly upper triangular matrices is itself strictly upper triangular, second-derivative operators such as $\mathbf{D}^{(d)}\mathbf{D}^{(r)}$ (used in the Laplacian for the Diffusion and Advection-Diffusion candidates) inherit the same structure, so the following conclusion holds independently of the order of the derivative operator involved.

A strictly upper triangular matrix of size $M$ satisfies $\mathbf{N}^{M}=\mathbf{0}$ \emph{exactly}.
This argument depends only on the strictly-upper-triangular shape of $\mathbf{N}$, requires neither an asymptotic limit nor any assumption on the magnitude of the matrix entries, and therefore applies unchanged to every candidate PDE considered in this study.

Consequently, every eigenvalue of $\mathbf{N}$ is exactly zero, and the equation-driven Koopman matrix $\mathbf{K}^{\star} = \exp(\Delta t\,\mathbf{N})$ has a single, degenerate eigenvalue $\lambda^{\star}=1$ with algebraic multiplicity $M$, \emph{independently of which candidate PDE generates $\mathbf{N}$}. Eigenvalue-based comparison between $\mathbf{K}^{\star}$ and $\hat{\mathbf{K}}$ is therefore uninformative for PDE identification: any observation-driven matrix $\hat{\mathbf{K}}$ whose eigenvalues cluster near unity appears equally consistent with every candidate model, irrespective of the true underlying dynamics. This motivates the use of a criterion that compares the action of the operators on data, such as the data-projected discrepancy, rather than their eigenvalues or eigenvectors alone.
\end{remark}

Fig.~\ref{fig:generator_structure} illustrates the nonzero pattern of $\mathbf{N}$ for a representative candidate (Advection-Diffusion, the only candidate with nonzero advection and diffusion coefficients simultaneously), for visual intuition only, and the exact nilpotency $\mathbf{N}^{M}=\mathbf{0}$ is already established for every candidate by Remark~\ref{rem:nilpotent}, and does not rely on this or any other numerical example.

Because every equation-driven matrix $\mathbf{K}^{\star}$ shares the same degenerate unit eigenvalue regardless of the underlying candidate PDE (Remark~\ref{rem:nilpotent}), neither the eigenvalues $\{\hat{\lambda}_j\}$, $\{\lambda_j^{\star}\}$ nor the corresponding Koopman modes 
\(
\hat{\boldsymbol{\xi}}_j
=\mathbf{C}^\intercal
\hat{\mathbf{r}}_j
\)
and
\(
\boldsymbol{\xi}_j^{\star}=
\mathbf{C}^\intercal\mathbf{r}_j^{\star}
\)
can discriminate between candidate models: any observation-driven matrix $\hat{\mathbf{K}}$ whose spectrum clusters near unity would appear equally consistent with every candidate, irrespective of the true governing dynamics. The proposed framework therefore compares the two operators through their action on the observed data rather than through their spectra.

Let
\begin{equation}
\hat{\mathbf{A}}_0
\coloneqq
\left(
\hat{\mathbf{a}}_0\;
\hat{\mathbf{a}}_1\;
\cdots\;
\hat{\mathbf{a}}_{N-2}
\right)
\end{equation}
denote the snapshot matrix of recovered Chebyshev coefficients introduced in Section~\ref{subsec:observation-driven}. The \emph{data-projected discrepancy} between the equation-driven matrix $\mathbf{K}^{\star}$ and the observation-driven matrix $\hat{\mathbf{K}}$ is defined as
\begin{equation}
d(\mathbf{K}^{\star},
\hat{\mathbf{K}})
\coloneqq
\frac{
\left\|
\left(
\mathbf{K}^{\star}
-
\hat{\mathbf{K}}
\right)
\hat{\mathbf{A}}_0
\right\|_{F}
}
{
\left\|
\mathbf{K}^{\star}
\hat{\mathbf{A}}_0
\right\|_{F}
},
\label{eq:data_projected_discrepancy}
\end{equation}
where $\|\cdot\|_{F}$ denotes the Frobenius norm. Unlike the eigenvalue- or eigenvector-based comparisons discussed above, the data-projected discrepancy directly evaluates the difference between the actions of the two Koopman matrices on the observed coefficient sequence rather than on their spectra, which are structurally degenerate for every candidate: a smaller value of $d$ indicates stronger agreement, while a larger value indicates a greater mismatch. This quantity remains well-defined and discriminative even though every candidate $\mathbf{K}^{\star}$ shares the same eigenvalue, and it is adopted as the sole criterion for PDE identification throughout the remainder of this paper.

\subsection{Observation Model Under Arbitrary Sampling}
\label{subsec:observation-model-under-arbitrary-sampling}

The formulations in Sections~\ref{subsec:observation-driven} to \ref{subsec:numerical-spectrum-linking} assume that the system state is represented in the Chebyshev spectral domain. In the ideal setting, observations are available directly at Chebyshev nodes, allowing the coefficient vector to be computed using DCT. In practical applications, however, measurements are generally acquired on observation grids that do not coincide with the Chebyshev nodes. Consequently, a mechanism is required to relate the Chebyshev spectral representation to observations obtained on arbitrary sampling grids.

Let $u \colon [-1,1]^D \to \mathbb{R}$ be a continuous distribution defined on the bounded domain
\([-1,1]^D.\)
In practice, only discrete observations of $u(\mathbf{q}), \mathbf{q}\in[-1,1]^D$ are available.
Using the Chebyshev basis introduced in Section~\ref{sec:theories}, the function can be represented as~\eqref{eq:cheb_finite_expansion}.

To derive the observation matrix, let
\(
\{\mathbf{q}_n\}_{n=0}^{N-1}
\)
denote the observation locations. Evaluating the Chebyshev expansion
\eqref{eq:cheb_finite_expansion}
at \(\mathbf{q}_n\) yields
\begin{equation}
v_n
=
u(\mathbf{q}_n)
\approx
\sum_{\mathbf{m}\in\mathcal{D}(\mathbf{M})}
a_{\mathbf{m}}
\tau_{\mathbf{m}}(\mathbf{q}_n).
\label{eq:pointwise_observation}
\end{equation}
Defining the Chebyshev basis evaluation matrix
\begin{equation}
[\boldsymbol{\Phi}]_{n,i} \coloneqq \tau_{\mathbf{m}_i}(\mathbf{q}_n),
\label{eq:chebysevevalmatrix}
\end{equation}
and collecting the pointwise observations into the vector
\(
\mathbf v = (v_0,v_1,\ldots,v_{N-1})^{\intercal}\in\mathbb{R}^N,
\)
the observation vector is expressed as
\begin{equation}
\mathbf{v} \approx \boldsymbol{\Phi}\mathbf{a}.
\label{eq:observation_phi_a}
\end{equation}

Let
\begin{equation}
\mathbf{a}
\coloneqq
(a_{\mathbf{m}})_{\mathbf{m}\in\mathcal{D}(\mathbf{M})}
\in
\mathbb{R}^{M},
\label{eq:coefficient_vector}
\end{equation}
denote the corresponding coefficient vector. Furthermore, let
\begin{equation}
\mathbf{x}
\coloneqq
(x_{\mathbf{m}})_{\mathbf{m}\in\mathcal{D}(\mathbf{M})}
\in
\mathbb{R}^{M},
\label{eq:chebyshev_sample_vector}
\end{equation}
denote the samples of the function on the Chebyshev nodes. The relationship between the coefficient vector $\mathbf{a}$ and the Chebyshev-node samples $\mathbf{x}$ is given by
\begin{equation}
\mathbf{x} \approx
\mathbf{C}^{\intercal}\mathbf{a},
\label{eq:chebyshev_node_samples}
\end{equation}
where $\mathbf{C}$ denotes the $D$-dimensional DCT-II matrix. \eqref{eq:chebyshev_node_samples} corresponds to the inverse DCT.

From \eqref{eq:chebyshev_node_samples} and the orthonormality of DCT-II, the coefficient vector $\mathbf{a}$ is written as
\(
\mathbf{a} \approx \mathbf{C}\mathbf{x}.
\)
Substituting this relation into
\eqref{eq:observation_phi_a}
gives
\(
\mathbf{v}
\approx
\boldsymbol{\Phi}\mathbf{C}\mathbf{x}.
\)
Since the observations $\mathbf{v}$ are expressed in the physical domain, whereas the unknown state $\mathbf{x}$ is represented by its samples on the Chebyshev nodes in~\eqref{eq:cheb_interior_nodes}, the DCT matrix naturally provides the mapping between these two representations.
Therefore, defining the observation matrix as
\begin{equation}
\mathbf{H}
\coloneqq
\boldsymbol{\Phi}\mathbf{C}
\label{eq:H_design_phiC}
\end{equation}
yields the observation model
\begin{equation}
\mathbf{v}
\approx
\mathbf{H}\mathbf{x},
\label{eq:observation_model_x}
\end{equation}
where
\(
\mathbf{H}
\in
\mathbb{R}^{N\times M},
\)
denotes the resampling geometry.

Substituting
\eqref{eq:chebyshev_node_samples}
into
\eqref{eq:observation_model_x}
yields
\begin{equation}
\mathbf{v} \approx
\mathbf{H}\mathbf{C}^{\intercal}\mathbf{a} = \boldsymbol{\Phi}\mathbf{a}
\label{eq:observation_model_a}
\end{equation}
\eqref{eq:observation_model_a} serves as the unified observation model throughout the remainder of this paper. Different sampling configurations are represented solely by different constructions of the Chebyshev basis evaluation matrix $\boldsymbol{\Phi}$, while the spectral representation and coefficient recovery formulation remain unchanged.

The specific structure of the Chebyshev basis evaluation matrix $\boldsymbol{\Phi}$ depends on the sampling configuration. In particular, Chebyshev, uniform, and irregular sampling schemes correspond to different constructions of $\boldsymbol{\Phi}$, as described in the following subsection.

\subsection{Observation Matrix Design}
\label{subsec:observation-matrix-design}

The observation matrix $\mathbf{H}$ introduced in \eqref{eq:observation_model_a} characterizes the relationship between the samples on the Chebyshev nodes and the measurements acquired on a given observation grid. In the proposed formulation, the observation process is represented by a linear mapping from the Chebyshev-node sample vector $\mathbf{x}$ to the observation vector $\mathbf{v}$.

Although the observation model introduced in Section~\ref{subsec:observation-model-under-arbitrary-sampling} is common to all sampling configurations, the construction of the observation matrix $\mathbf{H}$ depends on the observation geometry. The following subsections describe the corresponding construction for Chebyshev, uniform, and irregular sampling.

\subsubsection{Chebyshev Sampling}

Consider the ideal case in which observations are available directly on the Chebyshev nodes. In this situation, the observation vector coincides with the Chebyshev-node sample vector, i.e.,
\(
\mathbf{v} = \mathbf{x}.
\)

Therefore, the observation matrix becomes
\begin{equation}
\mathbf{H}
=
\mathbf{I}_{M},
\label{eq:H_chebyshev}
\end{equation}
where $\mathbf{I}_{M}$ denotes the $M\times M$ identity matrix.

Substituting \eqref{eq:H_chebyshev} into \eqref{eq:observation_model_a} gives
\(
\mathbf{v}
\approx
\mathbf{C}^{\intercal}\mathbf{a}.
\)
Hence, the original direct Chebyshev sampling formulation is recovered as a special case of the proposed unified observation model.

\subsubsection{Uniform Sampling}

In practical applications, measurements are often obtained on a uniform Cartesian grid rather than on the Chebyshev nodes. Measurements are acquired at the observation locations
\(
\{\mathbf{q}_n\}_{n=0}^{N-1}
\subset[-1,1]^D,
\)
which are arranged on a uniform Cartesian grid.

In this case, the observation matrix
\(\mathbf{H}\)
is constructed by evaluating the Chebyshev representation at the observation locations, thereby establishing the mapping from the Chebyshev-node samples to the uniformly sampled observations. Consequently, the observation model \eqref{eq:observation_model_x} provides a linear relationship between the Chebyshev-node samples and the measurements on the uniform grid.

The specific structure of $\boldsymbol{\Phi}$ depends on the sampling geometry and observation locations, while the underlying spectral representation remains unchanged.

\subsubsection{Irregular Sampling}

More generally, the observation locations
\(
\{\mathbf{q}_n\}_{n=0}^{N-1}
\subset[-1,1]^D
\)
may be distributed arbitrarily.

The same observation model can be applied by constructing a Chebyshev basis evaluation matrix $\boldsymbol{\Phi}$
that reflects the irregular sampling geometry. Compared with the uniform sampling case, only the observation locations differ, while the observation model and coefficient recovery formulation remain unchanged.

\begin{figure}[tb]
  \centering
  \includegraphics[width=.85\linewidth]{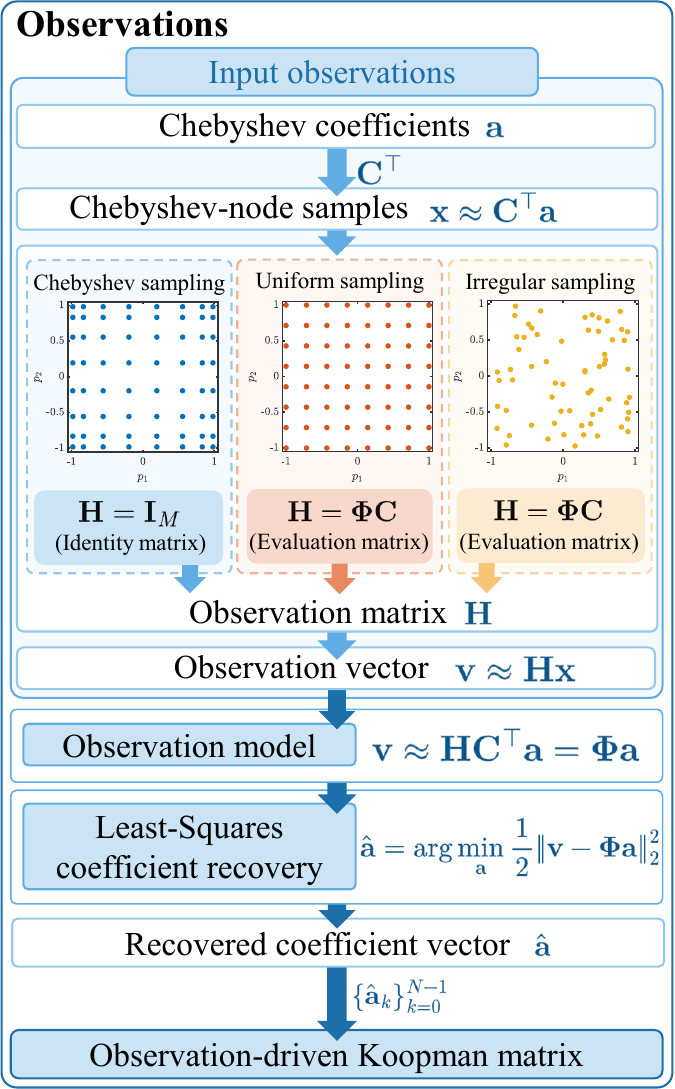}
  \caption{Proposed observation and coefficient recovery framework under arbitrary sampling. Observations acquired on Chebyshev, uniform, or irregular sampling grids are represented through the unified observation model $\mathbf{v} \approx \mathbf{H}\mathbf{C}^{\intercal}\mathbf{a} = \boldsymbol{\Phi}\mathbf{a}$, where the observation matrix $\mathbf{H}$ depends on the sampling configuration. The Chebyshev coefficient vector $\hat{\mathbf{a}}$ is recovered by solving the least-squares problem, and the resulting sequence of recovered coefficient vectors is subsequently used to estimate the observation-driven Koopman matrix.}
  
  \label{fig:observation_framework}
\end{figure}

Therefore, Chebyshev, uniform, and irregular sampling schemes can all be represented within the same mathematical framework. Fig.~\ref{fig:observation_framework} summarizes the unified observation framework considered in this work.

Once the Chebyshev basis evaluation matrix $\boldsymbol{\Phi}$ in~\eqref{eq:chebysevevalmatrix} has been specified, the problem of estimating the Chebyshev coefficient vector $\hat{\mathbf{a}}$ from the observations can be formulated as an inverse problem. This coefficient recovery problem is described in the following subsection.

\subsection{Coefficient Recovery}
\label{subsec:coefficient-recovery}

The Chebyshev coefficient vector $\hat{\mathbf{a}}$ can be estimated from the observations. 
Using the observation model in \eqref{eq:observation_model_a}, the coefficient recovery problem is formulated as the least-squares problem
\begin{equation}
\hat{\mathbf{a}}
\in
\arg\min_{\mathbf{a}}
\frac{1}{2}
\left\|
\mathbf{v}
-
\boldsymbol{\Phi}\mathbf{a}
\right\|_2^2.
\label{eq:coefficient_recovery_problem}
\end{equation}
This formulation seeks the coefficient vector whose predicted observations best match the measured observations under the observation model.
This condition on $N$ is therefore not an independent assumption but a direct consequence of the shape of $\boldsymbol{\Phi}$.

\eqref{eq:coefficient_recovery_problem} admits the closed-form expression
\begin{equation}
\hat{\mathbf{a}} = \boldsymbol{\Phi}^{\dagger}\mathbf{v},
\label{eq:coefficient_recovery_pseudoinverse}
\end{equation}
where 
\(\boldsymbol{\Phi}^\dagger\) denotes the Moore--Penrose pseudoinverse of \(\boldsymbol{\Phi}\), which reduces to \((\boldsymbol{\Phi}^\intercal \boldsymbol{\Phi})^{-1}\boldsymbol{\Phi}^\intercal\) when \(\boldsymbol{\Phi}\) has full column rank, to \(\boldsymbol{\Phi}^\intercal(\boldsymbol{\Phi}\boldsymbol{\Phi}^\intercal)^{-1}\) when \(\boldsymbol{\Phi}\) has full row rank, and, in the general (rank-deficient) case, to the non-unique \(\mathbf{V}\boldsymbol{\Sigma}^\dagger\mathbf{U}^\intercal\) via the singular value decomposition (SVD) \(\boldsymbol{\Phi}=\mathbf{U}\boldsymbol{\Sigma}\mathbf{V}^\intercal\)~\cite{Brunton2022Data-DrivenControl,Bramburger2024Data-DrivenSystems,Strang2019LinearData}.
This non-uniqueness is the source of the reduced robustness and occasional misidentification reported in Section~\ref{sec:performance}.
In either case, the estimated coefficient vector $\hat{\mathbf{a}}$ is subsequently used as the Chebyshev spectral representation of the observed state.

The full column rank requires $\boldsymbol{\Phi}$, and hence $\mathbf{HC}^{\intercal}$, to actually attain rank $M$, and $N\geq M$ alone does not guarantee this, as it can still fail for a poorly conditioned or degenerate sampling geometry. The dependence of this rank on the number and geometry of the observations is examined numerically in Section~\ref{subsec:observation-density}.

The proposed formulation naturally includes the direct Chebyshev sampling case, which corresponds to
\eqref{eq:H_chebyshev}
as discussed in Section~\ref{subsec:observation-matrix-design}.

Using the sequence of recovered coefficient vectors
\(
\{\hat{\mathbf{a}}_k\}_{k=0}^{N-1},
\)
the observation-driven Koopman matrix can be constructed by the least-squares procedure described in Section~\ref{subsec:observation-driven}. 
Consequently, the proposed formulation enables Koopman operator estimation from observations acquired on arbitrary sampling grids through the unified observation model and coefficient recovery process.

\section{Performance Evaluation}
\label{sec:performance}

In this section, the proposed numerical spectrum linking framework is evaluated under both ideal and non-ideal sampling conditions. The objective is to investigate whether the observation-driven Koopman operator $\hat{\mathbf{K}}$ remains consistent with the equation-driven Koopman operator $\mathbf{K}^{\star}$ when the governing dynamics are observed through different sampling configurations, including direct Chebyshev sampling, uniform sampling, and irregular sampling.

The experimental evaluation focuses on the effect of the sampling-dependent observation matrix on coefficient recovery, Koopman operator estimation, and PDE identification performance under different observation geometries.

The experiments are organized as follows.
Section~\ref{subsec:experimental-setup} describes the governing equations, numerical settings, and observation configurations.
Section~\ref{subsec:evaluation-metric} introduces the evaluation metric used for operator comparison.
Sections~\ref{subsec:chebyshev-sampling} to~\ref{subsec:irregular-sampling} present the results obtained under direct Chebyshev, uniform, and irregular sampling configurations, respectively.
Finally, Section~\ref{subsec:observation-density} investigates the influence of the observation density on coefficient recovery and PDE identification performance.

\subsection{Experimental Setup}
\label{subsec:experimental-setup}

To evaluate the proposed numerical spectrum linking framework under arbitrary sampling conditions, numerical experiments are conducted using a set of canonical PDEs representing distinct classes of spatiotemporal dynamics. Four governing equations of the form $\dot{u} = \mathcal{N}[u]$ are considered:
\begin{align}
\text{Advection-X:} \quad & \frac{\partial u}{\partial t} = a_x \frac{\partial u}{\partial x}, \\
\text{Advection-Y:} \quad & \frac{\partial u}{\partial t} = a_y \frac{\partial u}{\partial y}, \\
\text{Diffusion:} \quad & \frac{\partial u}{\partial t} = \nu \left(\frac{\partial^2 u}{\partial x^2} + \frac{\partial^2 u}{\partial y^2}\right), \\
\text{Advection-Diffusion:} \quad & \frac{\partial u}{\partial t} = a_x \frac{\partial u}{\partial x} + a_y \frac{\partial u}{\partial y} \nonumber\\
& \qquad + \nu \left(\frac{\partial^2 u}{\partial x^2} + \frac{\partial^2 u}{\partial y^2}\right),
\end{align}
where $a_x$ and $a_y$ denote the advection speeds in the $x$- and $y$-directions, respectively, and $\nu$ denotes the diffusivity. These equations serve as candidate models for the equation-driven Koopman operator $\mathbf{K}^{\star}$ construction described in Section~\ref{subsec:equation-driven}.

The spatial state is represented using a two-dimensional Chebyshev spectral expansion with truncation order $M_1=M_2=8$ in each spatial dimension, i.e., $\mathbf{M}=\operatorname{diag}(8,8)$, resulting in a spectral state dimension of $M = |\det(\mathbf{M})| = M_1 M_2 = 64$. This truncation order is fixed throughout the study. Temporal evolution is simulated using a sampling interval of $\Delta t = 5.0\times10^{-4}$ over the time horizon $T_{\mathrm{final}}=0.5$, corresponding to 1000 snapshots. To generate accurate reference trajectories, a finer internal integration step of $1.0\times10^{-5}$ is employed during numerical simulation.

\begin{figure}[tb]
  \centering
  \begin{minipage}[b]{0.459\textwidth}
    \centering
    \includegraphics[trim=10pt 30pt 10pt 5pt, clip, width=\linewidth]{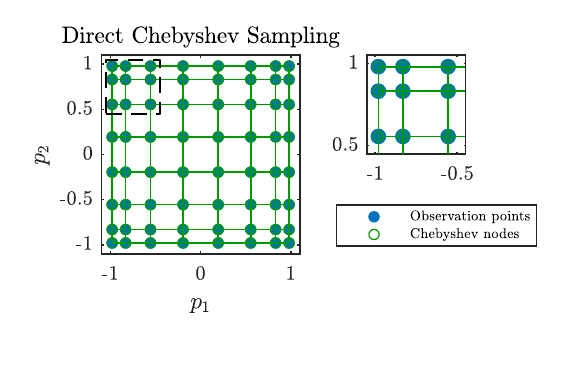}
    \label{fig:sampling_cheb}
    \centerline{\small (a) Direct Chebyshev sampling.}
  \end{minipage}
  \hfill
  \begin{minipage}[b]{0.459\textwidth}
    \centering
    \includegraphics[trim=10pt 30pt 10pt 5pt, clip, width=\linewidth]{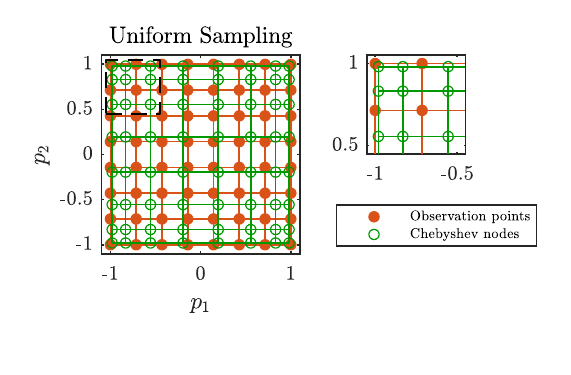}
    \label{fig:sampling_uniform}
    \centerline{\small (b) Uniform sampling.}
  \end{minipage}
  \hfill
  \begin{minipage}[b]{0.459\textwidth}
    \centering
    \includegraphics[trim=10pt 30pt 10pt 5pt, clip, width=\linewidth]{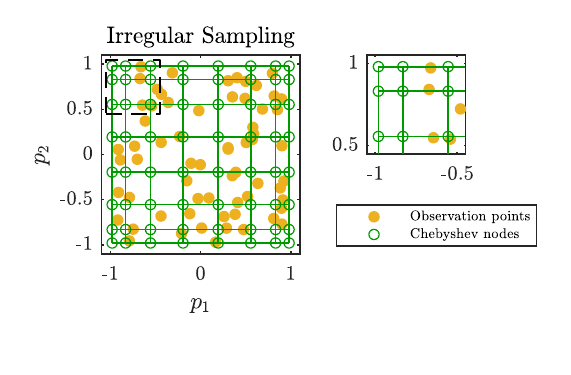}
    \label{fig:sampling_irregular}
    \centerline{\small (c) Irregular sampling.}
  \end{minipage}
  \caption{
    Observation geometries considered in this study. (a) Direct Chebyshev sampling: observation locations coincide with the Chebyshev nodes, resulting in no sampling mismatch, and the observation matrix $\mathbf{H}$ reduces to the identity matrix. (b) Uniform sampling: observations are acquired on a structured Cartesian grid that differs from the Chebyshev-node distribution. (c) Irregular sampling: observations are acquired at randomly distributed locations over the spatial domain. In (b) and (c), the resulting sampling mismatch is represented through the observation matrix $\mathbf{H}$, and the Chebyshev coefficients are recovered using the proposed least-squares formulation.
    }
  \label{fig:sampling_all}
\end{figure}

To investigate the effect of observation geometry on Koopman operator estimation and PDE identification, three observation configurations are considered: (a) direct Chebyshev sampling, (b) uniform sampling, and (c) irregular sampling. The corresponding observation geometries are illustrated in Fig.~\ref{fig:sampling_all}.

In Fig.~\ref{fig:sampling_all} (a) the direct Chebyshev sampling configuration, observations $\mathbf{v}$ are acquired directly on the Chebyshev nodes. Consequently, no sampling mismatch is present, and the observation matrix $\mathbf{H}$ reduces to the identity matrix, as defined in \eqref{eq:H_chebyshev}. This configuration serves as the ideal reference case and corresponds to the sampling assumption.

In Fig.~\ref{fig:sampling_all} (b) the uniform sampling configuration, observations $\mathbf{v}$ are collected on a structured Cartesian grid that differs from the Chebyshev-node distribution. The resulting sampling mismatch is represented through the observation matrix $\mathbf{H}$, and the Chebyshev coefficient vector $\hat{\mathbf{a}}$ is estimated through the coefficient recovery problem described in Section~\ref{subsec:coefficient-recovery}.

In Fig.~\ref{fig:sampling_all} (c) the irregular sampling configuration, observations $\mathbf{v}$ are acquired at randomly distributed locations over the spatial domain. Compared with the uniform sampling configuration, the irregular observation geometry generally produces a more challenging coefficient recovery problem because of the nonuniform distribution of the observation locations. As a result, the corresponding observation matrix $\mathbf{H}$ may become less well-conditioned, increasing the difficulty of accurate coefficient recovery.

The Chebyshev coefficient vectors $\hat{\mathbf{a}}$ are recovered using the least-squares formulation described in Section~\ref{subsec:coefficient-recovery} for all non-Chebyshev observation configurations. Unless otherwise stated, the number of observations is fixed at $N=64$, matching the spectral state dimension $M$. To further investigate the influence of observation density, the number of observations $N$ is varied in Section~\ref{subsec:observation-density}, while all other numerical settings and the least-squares coefficient recovery formulation remain unchanged.

\subsection{Evaluation Metric}
\label{subsec:evaluation-metric}

The consistency between the observation-driven Koopman matrix $\hat{\mathbf{K}}$ and the equation-driven Koopman matrix $\mathbf{K}^{\star}$ is evaluated using the data-projected discrepancy $d$ introduced in \eqref{eq:data_projected_discrepancy} for PDE identification throughout this paper.

For each observation configuration, the data-projected discrepancy is computed between the observation-driven Koopman matrix $\hat{\mathbf{K}}$ and the equation-driven Koopman matrices $\mathbf{K}^{\star}$ corresponding to all candidate PDEs. The candidate model yielding the minimum discrepancy is regarded as the identified governing equation. To further evaluate robustness, we use the identification margin, whose mathematical definition is given in Section~\ref{subsec:observation-density}: a positive margin indicates that the true governing PDE is uniquely identified, while a larger margin implies clearer separation from the competing candidate models.

\subsection{Direct Chebyshev Sampling}
\label{subsec:chebyshev-sampling}

We first consider the ideal observation setting in which measurements are acquired directly on the Chebyshev nodes. The least-squares coefficient recovery therefore reduces to the direct DCT-based computation. The top block of Table~\ref{tab:combined_discrepancy} reports the resulting data-projected discrepancy. As expected, the minimum discrepancy in each column is attained by the true governing PDE, with diagonal entries several orders of magnitude smaller than the off-diagonal entries, and all governing PDEs are therefore correctly identified. This configuration serves as the reference against which the non-Chebyshev sampling cases below are compared.

\subsection{Uniform Sampling}
\label{subsec:uniform-sampling}

We next consider the practical scenario in which observations are acquired on a uniform Cartesian grid rather than directly on the Chebyshev nodes, introducing a sampling mismatch between the observation grid and the spectral representation, and the Chebyshev coefficients are recovered via the least-squares formulation in \eqref{eq:coefficient_recovery_problem}. The middle block of Table~\ref{tab:combined_discrepancy} reports the resulting discrepancy, whose diagonal entries remain nearly identical to those under direct Chebyshev sampling. Consequently, all governing PDEs are correctly identified, demonstrating that the proposed least-squares recovery fully compensates for the uniform sampling mismatch.

\subsection{Irregular Sampling}
\label{subsec:irregular-sampling}

We finally consider the most challenging observation scenario, in which measurements are acquired at irregularly distributed locations, and since the observation geometry no longer preserves any structured relationship with the Chebyshev-node distribution, the coefficient recovery problem becomes more sensitive to the observation geometry and the resulting observation matrix may become less well-conditioned. The bottom block of Table~\ref{tab:combined_discrepancy} reports the discrepancy obtained from a single representative random draw at $N=64$ observations, and this choice is supported by the observation-density study of Section~\ref{subsec:observation-density}, which shows that at this observation count the identification margin is consistently positive with low seed-to-seed variance across 5 independent random draws (Table~\ref{tab:density_summary}), in contrast to the pronounced variability observed at lower densities. Although the diagonal entries are noticeably larger than under direct Chebyshev or uniform sampling, they remain substantially smaller than the off-diagonal entries, so all governing PDEs are still correctly identified. These results confirm that the proposed least-squares recovery remains effective even under nonuniform observation geometries.

\begin{table*}[tb]
    \centering
    \caption{Confusion matrix of the data-projected discrepancy $d$ (lower $\downarrow$ is better) across sampling configurations. Rows correspond to candidate PDEs and columns correspond to the true governing PDEs. The minimum discrepancy within each sampling block is highlighted in bold.}
    \label{tab:combined_discrepancy}
    \begin{tabular}{llcccc}
        \toprule
        Sampling & Candidate$\backslash$True & Advection-X & Advection-Y & Diffusion & Advection-Diffusion \\
        \midrule
        \multirow{4}{*}{Direct Chebyshev}
            & Advection-X   & \textbf{1.1646e-12} & 0.0083473 & 0.0094913 & 0.0031506 \\
            & Advection-Y   & 0.0055968 & \textbf{5.5232e-12} & 0.0059856 & 0.0035902 \\
            & Diffusion     & 0.0036225 & 0.0031936 & \textbf{4.0627e-07} & 0.0059185 \\
            & Advection-Diffusion & 0.0036305 & 0.0068674 & 0.0107582 & \textbf{4.0138e-07} \\
        \midrule
        \multirow{4}{*}{Uniform}
            & Advection-X   & \textbf{1.2467e-12} & 0.0083473 & 0.0094913 & 0.0031506 \\
            & Advection-Y   & 0.0055968 & \textbf{4.3935e-12} & 0.0059856 & 0.0035902 \\
            & Diffusion     & 0.0036225 & 0.0031936 & \textbf{3.7765e-07} & 0.0059185 \\
            & Advection-Diffusion & 0.0036305 & 0.0068674 & 0.0107582 & \textbf{2.8419e-07} \\
        \midrule
        \multirow{4}{*}{Irregular}
            & Advection-X   & \textbf{8.3096e-06} & 0.0083519 & 0.0094923 & 0.0031503 \\
            & Advection-Y   & 0.0055958 & \textbf{1.8841e-05} & 0.0059831 & 0.0035899 \\
            & Diffusion     & 0.0036213 & 0.0031925 & \textbf{2.0308e-05} & 0.0059182 \\
            & Advection-Diffusion & 0.0036300 & 0.0068723 & 0.0107582 & \textbf{4.5917e-07} \\
        \bottomrule
    \end{tabular}
\end{table*}

\subsection{Effect of Observation Density}
\label{subsec:observation-density}

The previous experiments demonstrate that the proposed framework correctly identifies the governing PDE under different sampling configurations when sufficient observation information is available. However, in practical applications, the number of available observations $N$ may be limited. Therefore, it is important to investigate how the observation density, parameterized by $N$, influences the robustness of coefficient recovery and PDE identification.

To this end, the number of observations $N$ is varied, while all other numerical parameters remain unchanged. For both uniform and irregular sampling configurations, $N$ is increased from $N=25$ to $N=121$. Because the irregular sampling locations are drawn at random, a single random draw may not be representative, and we therefore repeat the irregular sampling experiment over $S=5$ independent random seeds at each value of $N$ and report both the mean and the worst-case (least favorable seed) statistics. The uniform sampling configuration is deterministic for a given $N$ and therefore uses a single trial. The proposed least-squares coefficient recovery is applied in all cases, and the resulting observation-driven Koopman matrix $\hat{\mathbf{K}}$ is compared with the equation-driven Koopman matrices $\mathbf{K}^{\star}$ using the data-projected discrepancy $d$ introduced in Section~\ref{subsec:evaluation-metric}. The robustness of PDE identification is further evaluated using the identification margin $d_{\mathrm{margin}}$, defined later in this subsection \eqref{eq:identification_margin}.

\begin{figure}[tb]
    \centering
    \includegraphics[trim=15pt 5pt 10pt 10pt, clip, width=0.98\linewidth]{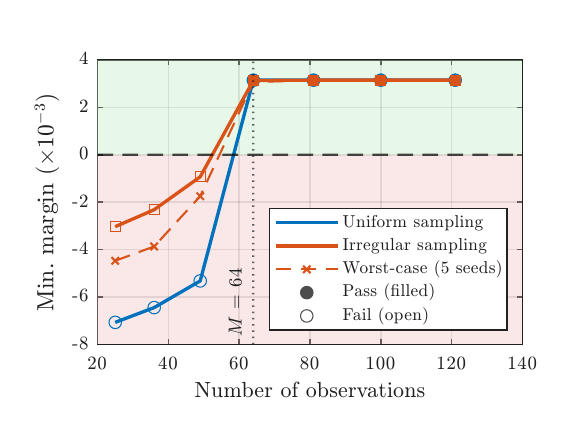}
    \caption{Identification margin $d_{\mathrm{margin}}$ as a function of the number of observations $N$ for uniform and irregular sampling configurations. For irregular sampling, the solid curve reports $d_{\mathrm{margin}}$ averaged over $S=5$ random seeds, and the dashed curve reports the worst-case margin among those seeds. A positive margin indicates that the projected discrepancy $d$ of the true governing PDE is smaller than that of every competing candidate. The transition at $N = M = 64$ observations corresponds to the spectral state dimension, where $\boldsymbol{\Phi}$ first attains full column rank, and notably, even the worst-case seed attains a positive margin at this threshold.}
    \label{fig:density_margin}
\end{figure}

Fig.~\ref{fig:density_margin} shows the identification margin $d_{\mathrm{margin}}$ as a function of the number of observations $N$. A positive margin indicates that the projected discrepancy $d$ of the true governing PDE is smaller than that of every competing candidate, whereas a negative margin indicates that at least one competing candidate yields a smaller projected discrepancy. For irregular sampling, the solid curve reports $d_{\mathrm{margin}}$ averaged over the $S=5$ random seeds, while the dashed curve reports the worst-case margin among those seeds.

The results in Fig.~\ref{fig:density_margin} and Table~\ref{tab:density_summary} reveal the following key observations:
\begin{itemize}
    \item \textbf{Below the spectral state dimension ($N < M = 64$):} $\boldsymbol{\Phi}$ is rank-deficient, resulting in reduced robustness and occasional misidentification of the governing PDE. Irregular sampling additionally exhibits pronounced seed-to-seed variability at these low densities, with the mean success rate $r$ rising from $r=35\%$ at $N=25$ to $r=70\%$ at $N=49$, and individual seeds ranging as widely as $50\%$--$100\%$ success at $N=49$ alone.
    \item \textbf{A clear transition occurs at $N=M=64$:} the matrix $\boldsymbol{\Phi}$ attains full column rank, enabling reliable recovery of the Chebyshev coefficients $\hat{\mathbf{a}}$, and $d_{\mathrm{margin}}$ becomes positive for both sampling configurations.
    \item \textbf{This transition holds even in the worst case:} the least favorable of the $S=5$ random seeds also achieves a positive margin once $N=M=64$, indicating that reliable identification at this threshold is not contingent on a fortunate random draw.
    \item \textbf{Beyond this threshold ($N>M$), performance plateaus:} increasing $N$ does not further increase the rank of $\boldsymbol{\Phi}$, since $\boldsymbol{\Phi}$ is an $N\times M$ matrix whose rank is mathematically bounded above by $\min(N,M)$, and once $N\ge M$, the rank ceiling of $M=64$ is guaranteed regardless of how many additional observations are collected, so $d_{\mathrm{margin}}$ remains nearly unchanged.
\end{itemize}
These observations indicate that reliable PDE identification depends primarily on obtaining a sufficient number of \emph{independent} observations rather than simply increasing the total number of measurements: once $\boldsymbol{\Phi}$ attains full column rank, additional observations beyond the spectral dimension $M$ provide only marginal improvement in identification robustness.

\begin{table*}[tb]
    \centering
    \caption{Summary of the effect of observation density on PDE identification. For irregular sampling, statistics are computed over 5 independent random seeds at each observation density, and the success rate is reported as the mean with the range (minimum--maximum) observed across seeds in parentheses, and uniform sampling is deterministic for a given number of observations and therefore uses a single trial. The rank of $\boldsymbol{\Phi}$ and the mean and worst-case identification margins are also reported for different observation densities.}
    \label{tab:density_summary}
    \begin{tabular}{lcccccc}
        \toprule
        Sampling & Number of Observations & Rank of $\boldsymbol{\Phi}$ & Success Rate & Mean Margin ($\times10^{-3}$) & Worst-case Margin ($\times10^{-3}$) & Trials \\
        \midrule
        \multirow{7}{*}{Uniform}
            & 25  & 25 & 25\%  & -7.07 & -7.07 & 1 \\
            & 36  & 36 & 25\%  & -6.45 & -6.45 & 1 \\
            & 49  & 49 & 25\%  & -5.32 & -5.32 & 1 \\
            & \textbf{64}  & \textbf{64} & \textbf{100\%} & \textbf{+3.15} & \textbf{+3.15} & 1 \\
            & 81  & 64 & 100\% & +3.15 & +3.15 & 1 \\
            & 100 & 64 & 100\% & +3.15 & +3.15 & 1 \\
            & 121 & 64 & 100\% & +3.15 & +3.15 & 1 \\
        \midrule
        \multirow{7}{*}{Irregular}
            & 25  & 25 & 35\% (25--50\%)  & -3.03 & -4.47 & 5 \\
            & 36  & 36 & 45\% (25--50\%)  & -2.32 & -3.87 & 5 \\
            & 49  & 49 & 70\% (50--100\%)  & -0.93 & -1.75 & 5 \\
            & \textbf{64}  & \textbf{64} & \textbf{100\% (100--100\%)} & \textbf{+3.13} & \textbf{+3.07} & 5 \\
            & 81  & 64 & 100\% (100--100\%) & +3.15 & +3.15 & 5 \\
            & 100 & 64 & 100\% (100--100\%) & +3.15 & +3.15 & 5 \\
            & 121 & 64 & 100\% (100--100\%) & +3.15 & +3.15 & 5 \\
        \bottomrule
    \end{tabular}
\end{table*}

The rank of $\boldsymbol{\Phi}$ indicates the number of independent constraints available for recovering the Chebyshev coefficient vector $\hat{\mathbf{a}}$. The success rate $r$ is defined as the percentage of governing PDEs that are correctly identified using the minimum data-projected discrepancy $d$. To further quantify the robustness of PDE identification, the identification margin $d_{\mathrm{margin}}$ is defined as
\begin{equation}
d_{\mathrm{margin}} = d_{\mathrm{comp}} - d_{\mathrm{true}},
\label{eq:identification_margin}
\end{equation}
where $d_{\mathrm{true}}$ denotes the projected discrepancy of the true governing PDE and $d_{\mathrm{comp}}$ denotes the smallest projected discrepancy among all competing candidate PDEs. Accordingly, a positive value of $d_{\mathrm{margin}}$ indicates successful identification of the true governing PDE. Table~\ref{tab:density_summary} summarizes these quantities for each observation density $N$.

\section{Conclusion}
\label{sec:conclusion}

A numerical spectrum linking framework is presented for PDE identification by comparing observation-driven and equation-driven Koopman matrices in a common Chebyshev spectral domain. To remove the restrictive assumption that observations are acquired directly on the Chebyshev nodes, a unified observation model together with a least-squares coefficient recovery formulation was introduced, enabling the proposed framework to operate with measurements obtained from arbitrary sampling configurations while remaining compatible with the direct Chebyshev sampling case.

The proposed framework was evaluated under direct Chebyshev, uniform, and irregular sampling configurations. The experimental results demonstrated that the proposed observation model and coefficient recovery formulation successfully recover the Chebyshev coefficient vectors $\hat{\mathbf{a}}$ and enable reliable PDE identification under both structured and unstructured observation geometries. Furthermore, the observation-density experiment showed that reliable PDE identification is consistently achieved once the matrix $\boldsymbol{\Phi}$ attains full column rank, providing a practical guideline for selecting the observation density required for robust coefficient recovery and Koopman operator estimation. In addition, we established that the equation-driven generator matrix $\mathbf{N}$ is structurally nilpotent for every candidate PDE considered, which forces the equation-driven Koopman matrix $\mathbf{K}^{\star}$ to a degenerate unit eigenvalue regardless of the underlying dynamics. This structural result motivates the adoption of the data-projected discrepancy, rather than purely eigenvalue-based criteria, as the identification criterion throughout this paper.

These results demonstrate that the proposed framework provides a unified and effective approach for Koopman-based PDE identification from arbitrarily sampled observations. By combining the observation model, least-squares coefficient recovery, and data-projected discrepancy within a common framework, the proposed method extends numerical spectrum linking to practical observation settings while preserving accurate identification performance.

Future work includes theoretical analysis of observation-induced perturbations in Koopman operator estimation, robustness evaluation under noisy and incomplete observations, replacement of the assumed Chebyshev basis with a data-driven basis transformation enabling PDE identification beyond a fixed candidate library, and extension of the proposed framework to higher-dimensional dynamical systems.

\bibliographystyle{IEEEtran}


\begin{thebibliography}{10}
\providecommand{\url}[1]{#1}
\csname url@samestyle\endcsname
\providecommand{\newblock}{\relax}
\providecommand{\bibinfo}[2]{#2}
\providecommand{\BIBentrySTDinterwordspacing}{\spaceskip=0pt\relax}
\providecommand{\BIBentryALTinterwordstretchfactor}{4}
\providecommand{\BIBentryALTinterwordspacing}{\spaceskip=\fontdimen2\font plus
\BIBentryALTinterwordstretchfactor\fontdimen3\font minus \fontdimen4\font\relax}
\providecommand{\BIBforeignlanguage}[2]{{%
\expandafter\ifx\csname l@#1\endcsname\relax
\typeout{** WARNING: IEEEtran.bst: No hyphenation pattern has been}%
\typeout{** loaded for the language `#1'. Using the pattern for}%
\typeout{** the default language instead.}%
\else
\language=\csname l@#1\endcsname
\fi
#2}}
\providecommand{\BIBdecl}{\relax}
\BIBdecl

\bibitem{Shannon1949a}
\BIBentryALTinterwordspacing
C.~E. Shannon, ``{Communication in the Presence of Noise},'' \emph{Proceedings of the Institute of Radio Engineers}, vol.~37, no.~1, pp. 10--21, 1949. [Online]. Available: \url{https://doi.org/10.1109/JRPROC.1949.232969}
\BIBentrySTDinterwordspacing

\bibitem{Unser2000a}
\BIBentryALTinterwordspacing
M.~Unser, ``{Sampling-50 years after Shannon},'' \emph{Proceedings of the IEEE}, vol.~88, no.~4, pp. 569--587, 2000. [Online]. Available: \url{https://doi.org/10.1109/5.843002}
\BIBentrySTDinterwordspacing

\bibitem{Eldar2014SamplingSystems}
\BIBentryALTinterwordspacing
Y.~C. Eldar, \emph{{Sampling Theory: Beyond Bandlimited Systems}}.\hskip 1em plus 0.5em minus 0.4em\relax Cambridge University Press, 12 2014. [Online]. Available: \url{https://doi.org/10.1017/CBO9780511762321}
\BIBentrySTDinterwordspacing

\bibitem{Zeng2024KoopmanSampling}
Z.~Zeng, J.~Liu, and Y.~Yuan, ``A generalized Nyquist-Shannon sampling theorem using the Koopman operator,'' \emph{IEEE Trans. Signal Process.}, vol.~72, pp.~3595--3610, 2024

\bibitem{Moteki2022CaptureBars}
\BIBentryALTinterwordspacing
D.~Moteki, T.~Murai, T.~Hoshino, H.~Yasuda, S.~Muramatsu, and K.~Hayasaka, ``{Capture method for digital twin of formation processes of sand bars},'' \emph{Physics of Fluids}, vol.~34, no.~3, 3 2022. [Online]. Available: \url{https://doi.org/10.1063/5.0085574}
\BIBentrySTDinterwordspacing

\bibitem{Moteki2023OnSandbars}
\BIBentryALTinterwordspacing
D.~Moteki, S.~Seki, S.~Muramatsu, K.~Hayasaka, and H.~Yasuda, ``{On the occurrence of sandbars},'' \emph{Physics of Fluids}, vol.~35, no.~1, 1 2023. [Online]. Available: \url{https://doi.org/10.1063/5.0128760}
\BIBentrySTDinterwordspacing

\bibitem{Ohara2024Physics-informedModel}
\BIBentryALTinterwordspacing
Y.~Ohara, D.~Moteki, S.~Muramatsu, K.~Hayasaka, and H.~Yasuda, ``{Physics-informed neural networks for inversion of river flow and geometry with shallow water model},'' \emph{Physics of Fluids}, vol.~36, no.~10, p. 106633, 10 2024. [Online]. Available: \url{https://doi.org/10.1063/5.0232852}
\BIBentrySTDinterwordspacing

\bibitem{Karisawa2025Acceleration-inducedSlopes}
\BIBentryALTinterwordspacing
H.~Karisawa and H.~Yasuda, ``{Acceleration-induced laminarization of sheet flow over smooth steep slopes},'' \emph{Physics of Fluids}, vol.~37, no.~8, p. 85198, 8 2025. [Online]. Available: \url{https://doi.org/10.1063/5.0279763}
\BIBentrySTDinterwordspacing

\bibitem{Brunton2022Data-DrivenControl}
\BIBentryALTinterwordspacing
S.~L. Brunton and J.~N. Kutz, \emph{{Data-Driven Science and Engineering: Machine Learning, Dynamical Systems, and Control}}.\hskip 1em plus 0.5em minus 0.4em\relax Cambridge University Press, 5 2022. [Online]. Available: \url{https://doi.org/10.1017/9781009089517}
\BIBentrySTDinterwordspacing

\bibitem{Bramburger2024Data-DrivenSystems}
\BIBentryALTinterwordspacing
J.~J. Bramburger, \emph{Data-Driven Methods for Dynamic Systems}. Philadelphia, PA: Society for Industrial and Applied Mathematics, 2024. [Online]. Available: \url{https://epubs.siam.org/doi/abs/10.1137/1.9781611978162}
\BIBentrySTDinterwordspacing
 
\bibitem{Strang2019LinearData}
\BIBentryALTinterwordspacing
G.~Strang, \emph{Linear Algebra and Learning from Data}. Wellesley, MA: Wellesley-Cambridge Press, 2019. [Online]. Available: \url{https://epubs.siam.org/doi/abs/10.1137/1.9780692196380}
\BIBentrySTDinterwordspacing

\bibitem{Baddoo2023Physics-informedDecomposition}
\BIBentryALTinterwordspacing
P.~J. Baddoo, B.~Herrmann, B.~J. McKeon, J.~Nathan~Kutz, and S.~L. Brunton, ``{Physics-informed dynamic mode decomposition},'' \emph{Proceedings of the Royal Society A}, vol. 479, no. 2271, p. 20220576, 3 2023. [Online]. Available: \url{https://doi.org/10.1098/rspa.2022.0576}
\BIBentrySTDinterwordspacing

\bibitem{Kobayashi2023Multi-ResolutionModeling}
E.~Kobayashi, H.~Yasuda, K.~Hayasaka, Y.~Otake, S.~Ono, and S.~Muramatsu, ``{Multi-Resolution Convolutional Dictionary Learning for Riverbed Dynamics Modeling},'' in \emph{ICASSP 2023 - 2023 IEEE International Conference on Acoustics, Speech and Signal Processing (ICASSP)}.\hskip 1em plus 0.5em minus 0.4em\relax Institute of Electrical and Electronics Engineers (IEEE), 5 2023, pp. 1--5.

\bibitem{2020TheControl}
\BIBentryALTinterwordspacing
A.~Mauroy, I.~Mezi{\'{c}}, and Y.~Susuki, Eds., \emph{{The Koopman Operator in Systems and Control}}, ser. Lecture Notes in Control and Information Sciences.\hskip 1em plus 0.5em minus 0.4em\relax Cham: Springer International Publishing, vol.~484, 2020. [Online]. Available: \url{https://doi.org/10.1007/978-3-030-35713-9}
\BIBentrySTDinterwordspacing

\bibitem{Williams2015ADecomposition}
\BIBentryALTinterwordspacing
M.~O. Williams, I.~G. Kevrekidis, and C.~W. Rowley, ``{A Data–Driven Approximation of the Koopman Operator: Extending Dynamic Mode Decomposition},'' \emph{Journal of Nonlinear Science}, vol.~25, no.~6, pp. 1307--1346, 12 2015. [Online]. Available: \url{https://doi.org/10.1007/s00332-015-9258-5}
\BIBentrySTDinterwordspacing

\bibitem{Korda2020Data-drivenOperator}
\BIBentryALTinterwordspacing
M.~Korda, M.~Putinar, and I.~Mezi{\'{c}}, ``{Data-driven spectral analysis of the Koopman operator},'' \emph{Applied and Computational Harmonic Analysis}, vol.~48, no.~2, pp. 599--629, 3 2020. [Online]. Available: \url{https://doi.org/10.1016/j.acha.2018.08.002}
\BIBentrySTDinterwordspacing

\bibitem{Hu2026HybridKoopman}
F.~Hu, H.-T.~Zhang, C.~Lv, and J.~Wang, ``Self-supervised Koopman operator learning for distributed final synchronization prediction of networked nonlinear dynamics,'' \emph{IEEE Trans. Neural Netw. Learn. Syst.}, vol.~37, no.~7, pp.~3325--3335, 2026

\bibitem{Wu2025DiscoveringFramework}
M.~Wu, W.~Li, L.~Yu, L.~Sun, J.~Liu, and W.~Li, ``{Discovering Mathematical Expressions Through DeepSymNet: A Classification-Based Symbolic Regression Framework},'' \emph{IEEE Transactions on Neural Networks and Learning Systems}, vol.~36, no.~1, pp. 1356--1370, 2025.

\bibitem{Lou2026Data-drivenSystems}
\BIBentryALTinterwordspacing
S.~Lou, H.~Xu, W.~Wang, L.~Lu, H.~Sun, Y.~Liu, L.~Zhang, D.~Zhang, and Y.~Chen,
``{Data-driven Discovery of Governing Differential Equations Across Physical Systems},''
6 2026. [Online]. Available:
\url{https://arxiv.org/abs/2606.09638}
\BIBentrySTDinterwordspacing

\bibitem{Rudy2017Data-drivenEquations}
\BIBentryALTinterwordspacing
S.~H. Rudy, S.~L. Brunton, J.~L. Proctor, and J.~N. Kutz, ``{Data-driven discovery of partial differential equations},'' \emph{Science Advances}, vol.~3, no.~4, 4 2017. [Online]. Available: \url{https://doi.org/10.1126/sciadv.1602614}
\BIBentrySTDinterwordspacing

\bibitem{Chen2022SymbolicSGA-PDE}
\BIBentryALTinterwordspacing
Y.~Chen, Y.~Luo, Q.~Liu, H.~Xu, and D.~Zhang, ``{Symbolic genetic algorithm for discovering open-form partial differential equations (SGA-PDE)},'' \emph{Physical Review Research}, vol.~4, no.~2, p. 023174, 6 2022. [Online]. Available: \url{https://doi.org/10.1103/PhysRevResearch.4.023174}
\BIBentrySTDinterwordspacing

\bibitem{Du2024DISCOVER:Learning}
\BIBentryALTinterwordspacing
M.~Du, Y.~Chen, and D.~Zhang, ``{DISCOVER: Deep identification of symbolically concise open-form partial differential equations via enhanced reinforcement learning},'' \emph{Physical Review Research}, vol.~6, no.~1, p. 013182, 1 2024. [Online]. Available: \url{https://doi.org/10.1103/PhysRevResearch.6.013182}
\BIBentrySTDinterwordspacing

\bibitem{Larranaga2026HowLow}
\BIBentryALTinterwordspacing
A.~Larra{\~n}aga, U.~Fasel, and S.~L.~Brunton, ``{How Low Can You Go? Active Learning for Sparse Model Discovery in the Ultra-Low-Data Limit},'' 6 2026. [Online]. Available: \url{https://arxiv.org/abs/2606.12182}
\BIBentrySTDinterwordspacing

\bibitem{Raissi2019Physics-informedEquations}
M.~Raissi, P.~Perdikaris, and G.~E. Karniadakis, ``{Physics-informed neural networks: A deep learning framework for solving forward and inverse problems involving nonlinear partial differential equations},'' \emph{Journal of Computational Physics}, vol. 378, pp. 686--707, 2 2019.

\bibitem{Huang2025PartialSurvey}
S.~Huang, W.~Feng, C.~Tang, Z.~He, C.~Yu, and J.~Lv, ``{Partial Differential Equations Meet Deep Neural Networks: A Survey},'' \emph{IEEE Transactions on Neural Networks and Learning Systems}, 2025.

\bibitem{Shi2024KoopmanStudyb}
\BIBentryALTinterwordspacing
D.~Shi and X.~Yang, ``{Koopman Spectral Linearization vs. Carleman Linearization: A Computational Comparison Study},'' \emph{Mathematics}, vol.~12, no.~14, p. 2156, 7 2024. [Online]. Available: \url{https://doi.org/10.3390/math12142156}
\BIBentrySTDinterwordspacing

\bibitem{Ahmed1974DiscreteTransform}
N.~Ahmed, T.~Natarajan, and K.~R. Rao, ``{Discrete Cosine Transform},'' \emph{IEEE Transactions on Computers}, vol. C-23, no.~1, pp. 90--93, 1974.

\bibitem{Ochoa-Dominguez2019DiscreteTransform}
H.~Ochoa-Dominguez and K.~R. Rao, \emph{{Discrete Cosine Transform}}, 2nd~ed.\hskip 1em plus 0.5em minus 0.4em\relax CRC Press, 2019.

\bibitem{Barrio2004AlgorithmsSeries}
\BIBentryALTinterwordspacing
R.~Barrio, ``{Algorithms for the integration and derivation of Chebyshev series},'' \emph{Applied Mathematics and Computation}, vol. 150, no.~3, pp. 707--717, 3 2004. [Online]. Available: \url{https://doi.org/10.1016/S0096-3003(03)00301-1}
\BIBentrySTDinterwordspacing

\bibitem{Sisaykeo2026NumSpecLink}
P.~Sisaykeo, and S.~Muramatsu, ``{Numerical Spectrum Linking: Identification of Governing PDE via Koopman-Chebyshev Approximation},'' in \emph{ICASSP 2026 - 2026 IEEE International Conference on Acoustics, Speech and Signal Processing (ICASSP)}.\hskip 1em plus 0.5em minus 0.4em\relax Institute of Electrical and Electronics Engineers (IEEE), 5 2026, pp. 1001--1005.

\bibitem{Strasser2025AnGuarantees}
\BIBentryALTinterwordspacing
R.~Str{\"a}sser, K.~Worthmann, I.~Mezi{\'{c}}, J.~Berberich, M.~Schaller, and F.~Allg{\"{o}}wer, ``{An overview of Koopman-based control: From error bounds to closed-loop guarantees},'' 9 2025. [Online]. Available: \url{https://arxiv.org/pdf/2509.02839}
\BIBentrySTDinterwordspacing

\bibitem{Bistrian2025ReducedLearning}
\BIBentryALTinterwordspacing
D.~A. Bistrian, ``{Reduced Order Data-Driven Twin Models for Nonlinear PDEs by Randomized Koopman Orthogonal Decomposition and Explainable Deep Learning},'' \emph{Mathematics}, vol.~13, no.~17, p. 2870, 9 2025. [Online]. Available: \url{https://doi.org/10.3390/math13172870}
\BIBentrySTDinterwordspacing

\bibitem{Mason2002ChebyshevPolynomials}
\BIBentryALTinterwordspacing
J.~Mason and D.~C. Handscomb, ``{Chebyshev Polynomials},'' 9 2002. [Online]. Available: \url{https://doi.org/10.1201/9781420036114}
\BIBentrySTDinterwordspacing

\bibitem{Karunakar2019ShiftedEquations}
\BIBentryALTinterwordspacing
P.~Karunakar and S.~Chakraverty, ``{Shifted Chebyshev polynomials based solution of partial differential equations},'' \emph{SN Applied Sciences}, vol.~1, no.~4, pp. 1--9, 4 2019. [Online]. Available: \url{https://doi.org/10.1007/s42452-019-0292-z}
\BIBentrySTDinterwordspacing

\bibitem{Prodinger2017RepresentingQuestions}
\BIBentryALTinterwordspacing
H.~Prodinger, ``{Representing derivatives of Chebyshev polynomials by Chebyshev polynomials and related questions},'' \emph{Open Mathematics}, vol.~15, no.~1, pp. 1156--1160, 1 2017. [Online]. Available: \url{https://doi.org/10.1515/math-2017-0096}
\BIBentrySTDinterwordspacing

\bibitem{Bedratyuk2022DerivationsPolynomials}
\BIBentryALTinterwordspacing
L.~P. Bedratyuk and N.~B. Lunio, ``{Derivations and Identities for Chebyshev Polynomials},'' \emph{Ukrainian Mathematical Journal}, vol.~73, no.~8, pp. 1175--1188, 1 2022. [Online]. Available: \url{https://doi.org/10.1007/s11253-022-01985-8}
\BIBentrySTDinterwordspacing

\bibitem{Ketthong2024Data-DrivenDisturbance}
\BIBentryALTinterwordspacing
P.~Ketthong, J.~Samkunta, N.~T. Mai, M.~A.~S. Kamal, I.~Murakami, and K.~Yamada, ``{Data-Driven Koopman Based System Identification for Partially Observed Dynamical Systems with Input and Disturbance},'' \emph{Sci}, vol.~6, no.~4, p.~84, 12 2024. [Online]. Available: \url{https://doi.org/10.3390/sci6040084}
\BIBentrySTDinterwordspacing

\end{thebibliography}

\vfill

\end{document}